\documentclass[reqno,12pt]{amsproc}
\usepackage{ucs}
\usepackage[english]{babel}
\usepackage{amssymb,amsmath}
\usepackage{indentfirst}
\usepackage[usenames]{color}
\usepackage{graphicx}
\usepackage[left=2cm,right=2cm,top=1.5cm,bottom=1.5cm,bindingoffset=0cm]{geometry}
\usepackage{amsthm}
\usepackage{amsfonts}
\usepackage{graphicx}
\graphicspath{{pictures/}}
\DeclareGraphicsExtensions{.pdf,.png,.jpg}
\usepackage{wrapfig}

\theoremstyle{plain}  
\newtheorem*{acknowledgements}{Acknowledgements}
\newtheorem{theorem}{\bf Theorem}
\newtheorem{lemma}{\bf Lemma}
\newtheorem{corollary}{\bf Corollary}
\newtheorem{prop}{\bf Proposition}

\newtheorem{remark}{\bf Remark}

\newcommand{\g}{\gamma}
\renewcommand{\a}{\alpha}

\usepackage{xcolor}
\usepackage{hyperref}

\usepackage{float}

\hypersetup{pdfstartview=FitH, linkcolor=linkcolor,urlcolor=urlcolor, colorlinks=true, linkcolor=blue, citecolor=blue}

\begin{document}
\title{Bounds for the number \\of multidimensional partitions}
\author{Kristina Oganesyan}
\address{Lomonosov Moscow State University, Moscow Center for fundamental and applied mathematics, Centre de Recerca Matem\`atica, Universitat Aut\`onoma de Barcelona}
\email{oganchris@gmail.com}
\thanks{This research was supported by the Russian Science Foundation (project no. 21-11-00131).}
\date{}

\begin{abstract}
We obtain estimates for the number $p_d(n)$ of $(d-1)$-dimensional integer partitions of a number $n$. It is known that the two-sided inequality $C_1(d)n^{1-1/d}<\log p_d(n)< C_2(d)n^{1-1/d}$ is always true and that $C_1(d)>1$ whenever $\log n> 3d$. However, establishing the $``$right$"$ dependence of $C_2$ on $d$ remained an open problem. We show that if $d$ is sufficiently small with respect to $n$, then $C_2$ does not depend on $d$, which means that $\log p_d(n)$ is up to an absolute constant equal to $n^{1-1/d}$. Besides, we provide estimates of $p_d(n)$ for different ranges of $d$ in terms of $n$, which give the asymptotics of $\log p_d(n)$ in each case.
\end{abstract}
\maketitle

\section{Introduction}
For a given $d$, we call a set $S\subset \mathbb{Z}_+^d$ {\it a lower set} (or {\it a downward closed set}) if for any $\mathbf{x}=(x_1,...,x_d)\in \mathbb{Z}_+^d$ the condition $\mathbf{x}\in S$ implies $\mathbf{x'}=(x'_1,...,x'_d)\in S$ for all $\mathbf{x'}\in \mathbb{Z}_+^d$ with $x'_i\leq x_i,\;1\leq i\leq d$. There is a one-to-one correspondence between $d$-dimensional lower sets of cardinality $n$ and $(d-1)$-dimensional partitions of $n$, that is, representations of the form 
$$n=\sum_{i_1=1}^{\infty}\sum_{i_2=1}^{\infty}...\sum_{i_{d-1}=1}^{\infty}n_{i_1i_2...i_{d-1}},\quad n_{i_1i_2...i_{d-1}}\in \mathbb{Z}_+,$$
where $n_{i_1i_2...i_{d-1}}\geq n_{j_1j_2...j_{d-1}}$ if $j_k\geq i_k$ for all $k=1,2,...,d-1.$ Thus, lower sets represent a geometric interpretation of integer partitions, or {\it partition diagrams}. By $p_d(n)$ we denote the number of lower sets in $\mathbb{Z}_+^d$ containing exactly $n$ points\footnote{In some sources the same value is denoted by $p_{d-1}(n)$.}.

\subsection{Small dimensional lower sets.} The history begins with finding the number $p_2(n)$ of integer partitions of a positive integer $n$, i.e. of representations of $n$ as a sum of nonincreasing positive integers, and evidently goes back to Leibniz \cite{L}. However, the first significant results in the partition theory were obtained much later by Euler \cite{E}. Any such representation $n=n_1+n_2+...n_k,\;n_1\geq n_2\geq ... \geq n_k$, or integer partition, can be visualized via the Young diagram, which consists of $n$ cells placed in $k$ rows and $n_1$ columns so that the $i$th row contains $n_i$ cells and the first cell in each row belongs to the first column. Another way to understand $p_2(n)$ is considering the following generating function \cite[Vol. 2, p.~1]{M2}
\begin{align*}
\prod_{k=1}^{\infty}(1-x^k)^{-1}=\sum_{n=0}^{\infty}p_2(n)x^n,
\end{align*}
where we assume $p_d(0)=1$ for any $d$. In 1917, Hardy and Ramanujan revealed the asymptotic behaviour of the function $p_2(n)$ (see \cite[(3)]{HR1} or \cite[(1.4)]{HR2}):
\begin{align*}
p_2(n)\sim \frac{e^{\sqrt{\frac{2n}{3}}\pi}}{4\sqrt{3}n}.
\end{align*}
Later, Rademacher \cite[(1.8)]{R} found an expansion of $p_2(n)$ as a convergent series.

In the case $d=3$, for the so-called plane partitions, the generating function was given by MacMahon 
 \cite{M1}:
\begin{align*}
\prod_{k=1}^{\infty}(1-x^k)^{-k}=\sum_{n=0}^{\infty}p_3(n)x^n
\end{align*}
(see \cite{Ch} for a simpler proof). The asymptotics of $p_3(n)$ was obtained by Wright \cite[(2.21)]{W}, namely,
\begin{align*}
p_3(n)\sim \frac{(2\zeta(3))^{\frac{7}{36}}e^{\zeta'(-1)}}{\sqrt{2\pi}n^{\frac{25}{36}}}e^{3(\zeta(3))^{\frac{1}{3}}2^{-\frac{2}{3}}n^{\frac{2}{3}}},
\end{align*}
where 
\begin{align*}
\zeta'(-1)=2\int\limits_0^{\infty}\frac{y \log y}{e^{2\pi y}-1}\;dy\approx -0.165421.
\end{align*}

For the cases $d>3$, no generating functions are known so far, although MacMahon conjectured that the function
\begin{align}\label{gen}
\prod_{k=1}^{\infty}(1-x^k)^{-\binom{d+n-2}{n-1}}
\end{align}
should generate $p_d(n)$ for every $d$, but this turned out to be wrong. On the other hand, some relations between the numbers $p_d(n)$ and the so-called MacMahon's numbers generated by \eqref{gen}, as well as some numerical values of $p_d(n)$, can be found in \cite{ABMM}. Besides, it was conjectured in \cite{MR} that MacMahon's numbers give the asymptotics of $\log p_d(n)/n^{1-1/d}$ for solid partitions, i.e. for $d=4$, and the hypothesis was accompanied by the exact values of $p_4(n),\;n\leq 50,$ and Monte Carlo simulations (see also \cite{BGP} for related numerical results in higher dimensions). However, the computations in \cite{DG} make this conjecture unlikely to be true for $d=4$. 

It is worth mentioning that an effective method for evaluating $p_4(n)$ is suggested in \cite{K}. Moreover, there is an algorithm that enables one to compute numbers of partitions for $n\leq 26$ in any dimension (see \cite{G}). 

Importantly, the partition theory has many applications in physics, as there are a lot of physical structures resembling that of multidimensional integer partitions. In particular, integer partitions are used to estimate the energy levels for a heavy nucleus \cite{BK} and to study the shape of crystal growth \cite{T}. Another direction of research is based on the existence of a one-to-one correspondence between partitions of an integer and microstates of a gas particles stored in a harmonic oscillator, not only in two-dimensional case \cite{AK,WH} but also in multidimensional setting~\cite{N}. 
 
Furthermore, certain classes of trigonometric polynomials associated with lower sets have recently turned out to be a powerful tool in multivariate approximation (see \cite{BDGJP, CMN, DPSTT} and references therein).

In the problem of estimating $p_d(n)$, the important relation
\begin{align}\label{bat}
C_1(d)\leq\frac{\log p_d(n)}{n^{1-\frac{1}{d}}}\leq C_2(d)
\end{align}
was established by Bhatia, Prasad and Arora \cite[(12), (16)]{BPA}, however the exact dependence of the constants on $d$ remained an open problem. Explicit values of $C_1$ and $C_2$ have recently been suggested in \cite[Th. 1.5]{DPSTT}, according to which \eqref{bat} holds with
\begin{align}\label{S}
C_1(d)=0.9\frac{d}{(d!)^{\frac{1}{d}}}\log 2,\quad C_2(d)=\pi\sqrt{\frac{2}{3}}d^{\log d},
\end{align}
where the upper bound holds for any $n\in\mathbb{N}$ and the lower bound is valid for $n>55^{d}$. Note that in this case $C_1(d)$ is uniformly bounded from below since Stirling's formula gives
\begin{align*}
d!<\sqrt{2\pi d}\Big(\frac{d}{e}\Big)^de^{\frac{1}{12d}}
\end{align*}
for all $d\geq 1$, and consequently, we have for $d\geq 3$,
$$C_1(d)\geq\frac{0.9e\log 2}{(2\pi de^{\frac{1}{6d}})^{\frac{1}{2d}}}\geq\frac{0.9e\log 2}{(6\pi e^{\frac{1}{18}})^{\frac{1}{6}}}>1.$$
So, for $n>55^d$, there holds $\log p_d(n)>n^{1-1/d}$ (see \cite[Section 2]{HR2} for the case $d=2$).

Another notable result (which appeared after the first version of our paper was written) in the recent preprint \cite{Y} states that
\begin{align*}
p_d(n)\leq d\zeta(d)^{\frac{1}{d}}n^{1-\frac{1}{d}}+(d-1)\log n
\end{align*}
holds for all $d\geq 2$ and $n$.

In the following theorem, we show that the constant $C_2$ is also independent of the dimension provided that $d$ is sufficiently small with respect to $n$.

\begin{theorem}\label{pr} For any $d\geq 2$ and any $n\geq (30d)^{2d^2}$, there holds
\begin{align*}
1<\frac{\log p_d(n)}{n^{1-\frac{1}{d}}}<7200.
\end{align*}
\end{theorem}

\subsection{High dimensional lower sets.} If we do not restrict ourselves to the case of a fixed (or relatively small) dimension $d$ and assume that $d$ grows somehow significantly along with $n$, then the general structure of lower sets changes, and the results of Theorem \ref{pr} are no longer true. Besides, estimate \eqref{bat} with $C_1$ and $C_2$ from \eqref{S} becomes quite rough if we just allow $d$ to be of order $\log n$. Somewhat better bounds for this setting were obtained in \cite[(24), (31)]{CMN}:
\begin{align*}
p_d(n)\leq 2^{dn}\quad\text{and}\quad p_d(n)\leq d^{n-1}(n-1)!
\end{align*}
for any positive integers $d$ and $n$. The latter inequality was strengthened and complemented by a lower bound in \cite[Th. 1.4]{DPSTT}
\begin{align*}
\binom{d+n-2}{n-1}\leq p_d(n)\leq d^{n-1}.
\end{align*}
Note that $\binom{d+n-2}{n-1}>d^{n-1}/(n-1)!$
\vspace{4pt}

The following theorem provides asymptotics of $\log p_d(n)$ for different orders of growth of $d$. 

\begin{theorem}\label{summing}
(a) If $d> n^3/2$, then
\begin{align*}
1\leq \frac{p_d(n)}{\binom{d+n-2}{d-1}}<\frac{1}{1-\frac{n^3}{2d}}.
\end{align*}

(b) If $dn^{-2}\to\infty$ as $n\to\infty$, then
\begin{align*}
\log p_d(n)=(n-1)(\log d-\log n+1)+o(n).
\end{align*}

(c) If $d$ satisfies $cn^2\leq d\leq Cn^2$ for some constants $c$ and $C$, then
\begin{align*}
\log p_d(n)=n\log n +O(n).
\end{align*}

(d) If $dn^{-2}\to 0$ and $\log d\geq \log n+o(\log n)$ as $n\to\infty$, then
\begin{align*}
\log p_d(n)=n\log n +o(n\log n).
\end{align*}
\end{theorem}
The paper is organized as follows. In Section \ref{sec1}, we prove Theorem \ref{pr} and some auxiliary statements that give an idea about the structure of lower sets in small dimensions. The case of high dimensions is treated in Section \ref{sec2}, where we deduce Theorem \ref{summing} from Propositions \ref{up}--\ref{4}, providing a more detailed analysis and explicit bounds.

\section{Lower sets in small dimensional spaces}\label{sec1}

\begin{figure}[H]
\center{\includegraphics[width=0.27\linewidth]{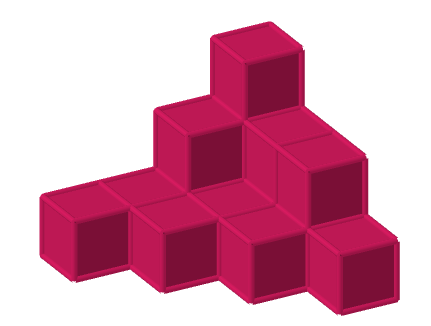}}
\caption{}
\label{pic}
\end{figure}

From now on we associate any point of a lower set with a unit cube having its center at that point. So, we will stick to this visualization of a lower set as a set of cubes leaning on one another. In Figure \ref{pic}, we give an example of such a visualization of a plane partition of $n=15$ with $$n_{11}=4,\;n_{12}=3,\;n_{13}=2,\;n_{14}=1,\;n_{21}=3,\;n_{22}=1,\;n_{31}=1,$$ so that the lower layer by itself represents a partition of $n_{11}+n_{12}+n_{13}+n_{14}=10$, while the next one, a partition of $n_{21}+n_{22}=4$.

For two cubes $q=(q_1,...,q_d)$ and $q'=(q'_1,...,q'_d)$, we write $q\succ q'$ if $q_i\geq q'_i$ for all $i=1,...,d$. If there holds either $q\succ q'$ or $q'\succ q$, we say that $q$ and $q'$ are {\it comparable}.

In this section we prove Theorem \ref{pr} and in the course of the proof reveal some features of the nature of lower sets. In the first place, we will be interested in the $``$top$"$ subsets of lower sets, which will play a crucial role in our further analysis. To be more specific, we need the following

{\bf Definition.} We call a subset $Q'$ of a lower set $Q$ {\it available} if for any $q'\in Q'$ there is no $q\in Q\setminus\{q'\}$ such that $q\succ q'$. 

In other words, it is such a subset that we can take its elements out in any order without breaking the lower set structure in any step. In terms of partially ordered sets, an available subset is an antichain of maximal elements of the lower set with respect to the partial order $\succ$. 

Denote by $M(Q)$ the maximal available subset of $Q$. The lemma below delivers the bound on $|M(Q)|$ that in some sense we would expect basing on our intuition: it seems that the more concentrated near the origin the lower set is, the richer its available subset can be (see Remark \ref{remark}). That is, if we construct a lower set $Q$ with too many cubes having too large sums of their coordinates, then we expect to have left unused many potential pairwise incomparable cubes $q\notin Q$ with smaller sums of coordinates. For instance, if a two-dimensional lower set $Q,\;|Q|=n,$ contains sufficiently many cubes with, say, the first coordinate considerably greater than $\sqrt{n}$, then there are a lot of columns having no available cube, so that we can enlarge $M(Q)$ by $``$descending$"$ the highly placed cubes and rearranging thereby the lower set $Q$.

\begin{lemma}\label{lem} For any $d\geq 2$ and $n\geq d^{6d\log d}$, for any $d$-dimensional lower set $Q,$ $|Q|=n,$ there holds
\begin{align*}
|M(Q)|\leq  \prod_{k=1}^{d-1}\bigg(1+\frac{1}{k^2}\bigg)n^{1-\frac{1}{d}}<\frac{\sinh \pi}{\pi}n^{1-\frac{1}{d}}.
\end{align*}
\end{lemma}

\begin{remark}\label{remark} Define the lower sets $S_k$ in $d$-dimensional space by the condition $x\in S_k\Leftrightarrow x_1+...+x_d\leq k$. Then it follows from Lemma \ref{lem} that the sets $S_k$, for large enough $k$, are optimal in the sense that their maximal available subsets are the largest possible up to an absolute constant. 
\end{remark}
 
\begin{proof}[Proof of Lemma \ref{lem}.] 
We prove the left-hand side inequality by induction on $d$, which will yield the assertion of the lemma. In the case $d=2$ we have $Q=Q_1\cup Q_2$, where $Q_i$ consists of all the cubes $q=(q_1,q_2)\in Q$ satisfying $q_i\leq \sqrt{n}-1,\;i=1,2$. Since $M(Q)$ cannot have more than one cube with a fixed $q_1$ or $q_2$, we can write $|M(Q)|\leq 2\sqrt{n}$. 

Suppose now that we proved the inequality for the dimensions $2,3,...,d-1,$ and let us prove it for $d\geq 3$. For simplicity we denote $$K_d:=\prod_{k=1}^{d-1}(1+k^{-2})<4$$ and divide the proof into several steps. 

{\bf Step 1. Obtaining a general bound for ${\bf |M(Q)|}$.} Consider all the nonempty subsets $Q_0,...,Q_m,\;m\leq n-1,$ being the intersections of $Q$ with the hyperplanes $q_1=0,...,m,$ respectively. They are lower sets themselves and if for some $j$ we have $q=(j,q_2,...,q_d)\in Q_j\cap M(Q)$, then $q=(j+s,q_2...,q_d)\notin Q_{j+s}$ for all $s\leq m-j$. Let $$n_i:=|Q_i|,\;0\leq i\leq m.$$ Note that we can apply the induction assumption to $Q_i$ with $n_i\geq (d-1)^{6(d-1)\log (d-1)}$. Now, taking into account that
\begin{align}\label{extra}
(d-1)^{6(d-1)\log (d-1)}< d^{6(d-1)(\log(d-1)-\log d)} n^{1-\frac{1}{d}}\leq d^{-3}n^{1-\frac{1}{d}},
\end{align} 
we have (assuming $n_{m+1}=0$)
\begin{align}\label{star}
|M(Q)|&\leq \sum_{i=0}^{m}\min\{n_i-n_{i+1},K_{d-1}n_i^{1-\frac{1}{d-1}}\}+\sum_{n_i\leq d^{-3}n^{1-\frac{1}{d}}}n_i-n_{i+1}\nonumber\\
&\leq \sum_{i=0}^{m}\min\{n_i-n_{i+1},K_{d-1}n_i^{1-\frac{1}{d-1}}\}+\frac{n^{1-\frac{1}{d}}}{d^3}\nonumber\\
&=:\sum_{i=0}^{m}\min\{\Delta_i, \Gamma_i\}+\frac{n^{1-\frac{1}{d}}}{d^3}=:\sum_{i=0}^{m}M_i+\frac{n^{1-\frac{1}{d}}}{d^3}=:F(n_0,...,n_m)+\frac{n^{1-\frac{1}{d}}}{d^3}.
\end{align}

If $m\leq n^{1/d}-1$, we have 
\begin{align}\label{n_m}
|M(Q)|\leq\max_{n_0+...n_m=n}\sum_{i=0}^{m}K_{d-1}n_i^{1-\frac{1}{d-1}}+\frac{n^{1-\frac{1}{d}}}{d^3}< K_{d}n^{1-\frac{1}{d}}
\end{align}
and there is nothing to prove. Thus, from now on, we can assume that 
\begin{align}\label{n_m_1}
m\geq n^{\frac{1}{d}}\qquad\text{and}\qquad n_m\leq n^{1-\frac{1}{d}}.
\end{align} 

We will maximize $F(n_0,...,n_m)$ over all $(n_0,...,n_m)$ in the set 
\begin{align*}
S_n:=\{(n_0,...,n_m)\in\mathbb{R}^+: n_0\geq ... \geq n_m,\; n_0+...+n_m=n\},
\end{align*}
permiting thereby $n_i$ to take noninteger values. We will see that an optimal $(n_0,...,n_m)$, roughly speaking, can have neither large values $n_i$ nor large differences $n_i-n_{i+1}$, which will in a way prescribe the behaviour of the sequence $n_i$. Take a point $(n_0,...,n_m)$ where this maximum is attained.

{\bf Step 2. Showing that ${\bf \Gamma_i\geq \Delta_i}$ for large ${\bf n_i}$.} First, by means of small perturbations of the optimal tuple $(n_0,...,n_m)$, we will show that while $n_i$ are sufficiently large, there must hold $\Delta_i\leq \Gamma_i$, so that the gaps between $n_i$ and $n_{i+1}$ cannot be too large in terms of $n_i$.

Assume that for some $i,\;0\leq i\leq m-1,$ such that $$n_{i}>8^{d-1},$$ we have $\Delta_i>\Gamma_i$. If we substitute the pair $(n_i,n_{i+1})$ by $(n_i-x,n_{i+1}+x)$ for sufficiently small positive $x$, the new point will still be in $S_n$ with $M_j,\;j\neq i-1,i,i+1,$ and $\Gamma_{i-1}$ remaining unchanged. At the same time, $\Delta_{i-1}$ increases, which means that $M_{i-1}$ does not decrease. Moreover, choosing $x$ small enough we can keep either the relation $\Delta_{i+1}\geq \Gamma_{i+1}$ or $\Delta_{i+1}\leq \Gamma_{i+1}$ true. Consider the two cases.

Case 1. $\Delta_{i+1}\geq \Gamma_{i+1}$. Then $M_i+M_{i+1}=\Gamma_i+\Gamma_{i+1}$ and the sum $\Gamma_i+\Gamma_{i+1}$ increases as $n_i$ and $n_{i+1}$ became closer to each other with their sum fixed. Thus, we increase $F(n_0,...,n_m)$, which contradicts the definition of $(n_0,...,n_m)$.

Case 2. $\Delta_{i+1}\leq \Gamma_{i+1}$. The value $\Delta_{i+1}+\Gamma_i$ changes in 
\begin{align*}
x-K_{d-1}\big(n_i^{1-\frac{1}{d-1}}-(n_{i}-x)^{1-\frac{1}{d-1}}\big)\geq x\Big(1-K_{d-1}\frac{d-2}{d-1}(n_{i}-x)^{-\frac{1}{d-1}}\Big)>0,
\end{align*}
since $n_i\geq 8^{d-1}>K_{d-1}^{d-1}$. Hence, $M_{i}+M_{i+1}$ increases. Thus, we increase $F(n_0,...,n_m)$, which once again contradicts the definition of $(n_0,...,n_m)$.

The fact that both cases led us to contradictions means that there holds 
\begin{align}\label{g}
\Delta_{i}\leq \Gamma_{i},\qquad 0\leq i\leq \min\{p,m-1\},
\end{align}
where $p$ is the maximal index satisfying $n_p\geq 8^{d-1}$.

The rest of the proof we divide into two cases: the case of $``$large$"$ and the case of $``$small$"$ values $n_0$. We will see that $n_0$ cannot be large at a point of the maximum of $F(n_0,...,n_m)$.

{\bf Step 3. Proving that ${\bf n_0}$ cannot be large.} Assume that $n_0> 2K_{d-1}n^{1-1/d}$. Define the sequence $\{a_i\}_{i=0}^{\infty}$ in the following way: 
\begin{align}\label{def_a}
a_0:=n_0\quad\text{ and}\quad a_{i}=a_{i-1}-K_{d-1}a_{i-1}^{1-\frac{1}{d-1}}\;\text{for}\;i\geq 1.
\end{align}
Let us estimate the maximal number $k$ such that
\begin{align*}
a_k>\frac{K_{d-1}n^{1-\frac{1}{d}}}{2}\quad\text{and}\quad \sum_{i=0}^{k} a_i\leq n.
\end{align*}
From the definition of $a_i$ we see that the ratio $a_i/a_{i+1}$ increases along with $i$. So,
\begin{align*}
\frac{a_{k/2}}{0.5K_{d-1}n^{1-\frac{1}{d}}}>\frac{a_{k/2}}{a_{k}}>\frac{a_0}{a_{k/2}}> \frac{2K_{d-1}n^{1-\frac{1}{d}}}{a_{k/2}},
\end{align*}
where in the case of odd $k$ we understand $a_{k/2}$ as $(a_{(k-1)/2}+a_{(k+1)/2})/2$. Hence, $a_{k/2}>K_{d-1}n^{1-1/d}.$ Since $a_i-a_{i+1}$ decreases, we have $(k+1)a_{k/2}\leq\sum_{i=0}^ka_i\leq n,$ which yields $$k+1<n^{\frac{1}{d}}K_{d-1}^{-1}.$$ So,
\begin{align*}
a_0-a_{k+1}=\sum_{i=0}^kK_{d-1}a_i^{1-\frac{1}{d-1}}\leq K_{d-1}\frac{n^{\frac{1}{d}}}{K_{d-1}}&\bigg(\frac{nK_{d-1}}{n^{1/d}}\bigg)^{1-\frac{1}{d-1}}=K_{d-1}^{1-\frac{1}{d-1}}n^{1-\frac{1}{d}}< K_{d-1}n^{1-\frac{1}{d}},
\end{align*}
whence $a_{k+1}>K_{d-1}n^{1-1/d}>0.5K_{d-1}n^{1-1/d}.$ Thus, the sum of $a_i$ becomes equal to $n$ before $a_i$ reaches $0.5K_{d-1}n^{1-1/d}$. Therefore, according to \eqref{g}, since $n_i$ for $i\leq p$ decreases slower than $a_i$ does, we obtain that $p=m$ and $$m+1\leq k+1<n^{\frac{1}{d}}K_{d-1}^{-1}<n^{\frac{1}{d}},$$ which contradicts \eqref{n_m_1}. Thus,
$$n_0\leq 2K_{d-1}n^{1-1/d}.$$

{\bf Step 4. Proving that ${\bf n_0-n_p}$ and ${\bf n_p}$ cannot be large.} Assume that 
\begin{align}\label{assume}
n_0-n_p>L_dn^{1-\frac{1}{d}},\quad L_d:=K_{d-1}\Big(1+\frac{2}{3(d-1)^{2}}\Big).
\end{align}
Considering the sequence $\{a_i\}$ given by \eqref{def_a}, denote by $q$ the maximal index such that $$a_{q}\geq 8^{d-1}\quad\text{and}\quad\sum_{i=0}^q a_i\leq n.$$ We divide the interval $(a_{q}+\varepsilon,a_0]$ into $$I_j:=(\nu_j,\mu_j]:=(A_jn^{1-\frac{1}{d}},(A_j+n^{-\frac{1}{2d}})n^{1-\frac{1}{d}}],\quad A_j=A_{j+1}+n^{-\frac{1}{2d}},$$ where $\varepsilon\in [0,n^{1-3/(2d)})$ is chosen uniquely. Note that $|I_j|=n^{1-3/(2d)}$ for all $j$. Denote the number of $a_i$'s belonging to $I_j$ by $k_{A_j}$ and let $a_{i_j}$ be the greatest of $a_i$ that belongs to $I_j$. 

Now, in order to prove that the assumption \eqref{assume} cannot hold, we are going to show that each $I_j$ contains significantly many terms $a_i$, and this will yield that $a_i$, and therefore $n_i$, cannot decrease considerably until the sum of its first terms becomes equal to $n$. The fact that $a_0=n_0$ is not very large provides some kind of regularity of $a_i$, namely, it will ensure that the leaps between $a_i$ are small enough to get appropriate estimates on $k_{A_j}$ for each of the intervals $I_j$.

Fix some $j$ and suppress for simplicity the index $j$ in $A_j$. Suppose that 
\begin{align}\label{sup}
k_A< \frac{A^{\frac{1}{d-1}-1}n^{\frac{1}{2d}}}{L_d}.
\end{align}
Since the ratio $a_i/a_{i+1}$ increases and $a_{i_{j+1}-1}>\nu_j=A_jn^{1-1/d}$, we have
\begin{align}\label{obs}
\bigg(\frac{1}{1-K_{d-1}A^{-\frac{1}{d-1}}n^{-\frac{1}{d}}}\bigg)^{k_A}&\geq \bigg(\frac{1}{1-K_{d-1}a_{i_{j+1}-1}^{-\frac{1}{d-1}}}\bigg)^{k_A}\geq \Big(\frac{a_{i_{j+1}-1}}{a_{i_{j+1}}}\Big)^{k_A} \geq \frac{a_{i_j}}{a_{i_{j+1}}}\geq \frac{\mu_j-K_{d-1}a_{i_j-1}^{1-\frac{1}{d-1}}}{\mu_{j+1}}\nonumber\\
&\geq \frac{(A+n^{-\frac{1}{2d}})n^{1-\frac{1}{d}}-K_{d-1}a_0^{1-\frac{1}{d-1}}}{An^{1-\frac{1}{d}}}\geq 1+\frac{n^{-\frac{1}{2d}}}{A}-\frac{K_{d-1}^22n^{-\frac{1}{d}}}{A}.
\end{align}
At the same time, since $(k_A+1+x)/(2+x)\leq (k_A+1)/2$ for $x\geq 0$, the ratio between consequent terms in the binomial expansion of 
\begin{align*}
(1-K_{d-1}A^{-\frac{1}{d-1}}n^{-\frac{1}{d}})^{-k_A}
\end{align*}
is less than
\begin{align*}
\frac{k_A+1}{2}K_{d-1}A^{-\frac{1}{d-1}}n^{-\frac{1}{d}}&<\frac{K_{d-1}A^{-\frac{1}{d-1}}n^{-\frac{1}{d}}}{2}+\frac{A^{-1}n^{-\frac{1}{2d}}K_{d-1}}{2L_d}\\
&\leq \frac{K_{d-1}(d-1)^{\frac{2}{d-1}}n^{-\frac{1}{d}}}{2}+\frac{(d-1)^{2}n^{-\frac{1}{2d}}K_{d-1}}{2L_d}\leq (d-1)^{2}n^{-\frac{1}{2d}}.
\end{align*}
This implies
\begin{align*}
(1-K_{d-1}A^{-\frac{1}{d-1}}n^{-\frac{1}{d}})^{-k_A}&< 1+k_AK_{d-1}A^{-\frac{1}{d-1}}n^{-\frac{1}{d}}\frac{1}{1-(d-1)^{2}n^{-\frac{1}{2d}}}\\
&<1+\frac{n^{-\frac{1}{2d}}}{A}\frac{1}{\Big(1+\frac{2}{3(d-1)^2}\Big)(1-(d-1)^{2}n^{-\frac{1}{2d}})}.
\end{align*}
Combining this with \eqref{obs}, we obtain
\begin{align*}
\frac{1}{\Big(1+\frac{2}{3(d-1)^2}\Big)(1-(d-1)^{2}n^{-\frac{1}{2d}})}>1-2n^{-\frac{1}{2d}}K_{d-1}^2,
\end{align*}
which yields
\begin{align*}
0&> -2n^{-\frac{1}{2d}}K_{d-1}^2+\frac{2}{3(d-1)^{2}}-\frac{4n^{-\frac{1}{2d}}K_{d-1}^2}{3(d-1)^2}-(d-1)^2n^{-\frac{1}{2d}}+2(d-1)^2n^{-\frac{1}{d}}K_{d-1}^2\\
&-\frac{2n^{-\frac{1}{2d}}}{3}+\frac{4n^{-\frac{1}{d}}K_{d-1}^2}{3}\geq \frac{2}{3(d-1)^2}-\frac{68}{3}n^{-\frac{1}{2d}}(d-1)^2\geq 0,
\end{align*}
since $n>d^{18d}>34^{2d}(d-1)^{8d}$. This contradiction disproves \eqref{sup}, whence
\begin{align*}
k_{A_j}\geq\frac{A_j^{\frac{1}{d-1}-1}n^{\frac{1}{2d}}}{L_d}.
\end{align*} 

Summing up this inequality over all $j$, we derive 
\begin{align*}
n\geq \sum_{a_i\in\cup I_j}a_i\geq \sum_{j}k_{A_j}A_jn^{1-\frac{1}{d}}\geq \sum_{j}\frac{n^{1-\frac{1}{2d}}A_j^{\frac{1}{d-1}}}{L_d}&> nL_d^{-1}\int\limits_{a_qn^{\frac{1}{d}-1}+n^{-\frac{1}{2d}}}^{a_0n^{\frac{1}{d}-1}-n^{-\frac{1}{2d}}}x^{\frac{1}{d-1}}\;dx\\
&=:nL_d^{-1}\int\limits_{y}^{y+z-2n^{-\frac{1}{2d}}}x^{\frac{1}{d-1}}\;dx,
\end{align*}
where $z\geq L_d$ by the assumption \eqref{assume}. This means that there holds
\begin{align}\label{cont}
1> L_d^{-1}\int\limits_{y}^{y+L_{d}-2n^{-\frac{1}{2d}}}x^{\frac{1}{d-1}}\;dx&\geq L_d^{-1}\frac{d-1}{d}(L_d-2n^{-\frac{1}{2d}})^{\frac{d}{d-1}}\nonumber\\
&\geq (1-2n^{-\frac{1}{2d}})L_d^{\frac{1}{d-1}}\frac{d-1}{d}\geq \frac{L_d^{\frac{1}{d-1}}\frac{d-1}{d}}{1+3n^{-\frac{1}{2d}}}\geq 1.
\end{align}
Let us prove the latter inequality. It suffices to show that 
\begin{align*}
L_d\geq \Big(1+\frac{1}{d-1}\Big)^{d-1}(1+3n^{-\frac{1}{2d}})^{d-1}.
\end{align*}
This, in turn, will follow from
\begin{align}\label{1}
1+\frac{2}{3(d-1)^2}\geq (1+3n^{-\frac{1}{2d}})^{d-1}\Big(1+\frac{1}{2(d-1)^2}\Big)
\end{align}
and 
\begin{align}\label{2}
K_{d-1}\Big(1+\frac{1}{2(d-1)^2}\Big)\geq \Big(1+\frac{1}{d-1}\Big)^{d-1}.
\end{align}
Firstly, by the assumption of the lemma we have $n> 6^{2d}(d-1)^{2d}$, so
\begin{align*}
1+\frac{2}{3(d-1)^2}&> 1+\frac{1}{2(d-1)^2}+6(d-1)n^{-\frac{1}{2d}}+\frac{3n^{-\frac{1}{2d}}}{d-1}\\
&=(1+6(d-1)n^{-\frac{1}{2d}})\Big(1+\frac{1}{2(d-1)^2}\Big)\geq (1+3n^{-\frac{1}{2d}})^{d-1}\Big(1+\frac{1}{2(d-1)^2}\Big),
\end{align*}
which proves \eqref{1}. Secondly, note that for $d=3$ both sides of \eqref{2} are equal $2.25$ and for $d=4$ inequality \eqref{2} becomes $95/36\geq (4/3)^3$. For $d\geq 5$, \eqref{2} follows from the fact that the left-hand side is greater than $e$. 

The contradiction in \eqref{cont} along with the fact that $a_i$'s decrease faster than $n_i$ shows that \eqref{assume} does not hold and therefore
\begin{align}\label{n0p}
n_0-n_p\leq L_dn^{1-\frac{1}{d}}.
\end{align}
In addition, we note that if $p\neq m$, then by \eqref{g}
\begin{align*}
n_{p+1}\geq n_p(1-K_{d-1}n_p^{-\frac{1}{d-1}})>n_p\Big(1-K_{d-1}(8^{d-1}))^{-\frac{1}{d-1}}\Big)>0.5n_p,
\end{align*}
whence $n_p<2\cdot 8^{d-1}$. 

{\bf Step 5. Final estimates.} Finally, in light of \eqref{star}, \eqref{n_m}, and \eqref{n0p}, we obtain
\begin{align*}
|M(Q)|&\leq \frac{n^{1-\frac{1}{d}}}{d^3}+\sum_{i=0}^{p-1}(n_i-n_{i+1})+\sum_{i=p}^{m}M_i\\
&< \frac{n^{1-\frac{1}{d}}}{d^3}+(n_0-n_p)+2\cdot 8^{d-1}+K_{d-1}n_m^{1-\frac{1}{d-1}}\nonumber\\
&< \frac{n^{1-\frac{1}{d}}}{d^3}+L_dn^{1-\frac{1}{d}}+2\cdot 8^{d-1}+K_{d-1}n^{1-\frac{2}{d}}\nonumber\\
&<K_{d-1}n^{1-\frac{1}{d}}\bigg(\frac{1}{2d^3}+\Big(1+\frac{2}{3(d-1)^2}\Big)+\frac{8^{d-1}}{n^{1-\frac{1}{d}}}+n^{-\frac{1}{d}}\bigg)\nonumber\\
&\leq K_{d-1}n^{1-\frac{1}{d}}\Big(1+\frac{2}{3(d-1)^2}+\frac{1}{2d^3}+d^{-10(d-1)\log d}+n^{-\frac{1}{d}}\Big)< K_dn^{1-\frac{1}{d}},
\end{align*}
since $n> 12^d(d-1)^{2d}$, and the proof of the lemma is complete.

In conclusion, we just remark that under the assumptions of the lemma the argument above gives
\begin{align}\label{xi}
\max_{(n_0,...,n_m)\in S_n}\sum_{i=0}^{m}\min\{n_i-n_{i+1},K_{d-1}n_i^{1-\frac{1}{d-1}}\}\leq K_{d-1}n^{1-\frac{1}{d}}\Big(1+\frac{2}{3(d-1)^2}+d^{-10(d-1)\log d}+n^{-\frac{1}{d}}\Big).
\end{align}
\end{proof}

\begin{remark} Without any restriction on $d$ and $n$ we can straightforwardly show that there always holds
\begin{align*}
|M(Q)|\leq dn^{1-\frac{1}{d}}.
\end{align*}
\end{remark}
 
\begin{proof}
Proceeding by induction as in the proof of Lemma \ref{lem} we obtain 
\begin{align*}
|M(Q)|&\leq \max_{n_0+...+n_m=n} \sum_{i=0}^{m}\min\{n_i-n_{i+1},(d-1)n_i^{1-\frac{1}{d-1}}\}.
\end{align*}
As long as $Q$ is a lower set, there holds $n_i\geq n_{i+1}$ for any $i=1,...,m-1,$ so $n_{\lfloor n^{1/d}\rfloor}\leq n^{1-1/d}.$ Thus, $$\sum_{k\geq \lfloor n^{1/d}\rfloor}(n_k-n_{k+1})\leq n^{1-\frac{1}{d}}$$ and
\begin{align*}
|M(Q)|\leq n^{1-\frac{1}{d}}+(d-1)\max_{n_0+...+n_{\lfloor n^{1/d}\rfloor -1}\leq n} \sum_{k=0}^{\lfloor n^{\frac{1}{d}}\rfloor-1}n_k^{1-\frac{1}{d-1}}\leq n^{1-\frac{1}{d}}+(d-1) n^{\frac{1}{d}} \Big(\frac{n}{n^{1/d}}\Big)^{1-\frac{1}{d-1}}= dn^{1-\frac{1}{d}}.
\end{align*}
\end{proof}

Now, as we already have the bound for the cardinalities of the available subsets, we are able to get some estimates for the number of lower subsets of a lower set. Define $$T(n):=\max_{\text{lower sets}\;Q:\;|Q|=n}|M(Q)|.$$

\begin{lemma}\label{cor} For the number $C(Q,k,d)$ of all lower subsets $Q',\;|Q'|\geq n-k$, of a lower set $Q,\;|Q|=n,$ in $d$-dimensional space
 there holds
\begin{align*}
C(Q,k,d)< \Big(\max\Big\{8,\frac{4eT(n)}{k}\Big\}\Big)^k.
\end{align*}
\end{lemma}

\begin{proof}
First we show that every lower subset $Q'$ of a lower set $Q$ can be constructed by successively discarding cubes of $Q$ one by one so that in any step the current set remains being a lower set. 

Indeed, let us list all the cubes we have to discard from $Q$ in a sequence in an arbitrary order. By a {\it disorder} we call a pair $(q,q')$ of cubes in this sequence such that $q$ goes after $q'$ in it, but $q\succ q'$. Now, if there is a disorder $(q,q')$ in the sequence, we just swap $q$ and $q'$, which removes the disorder and does not create any new one. Thus, we can rearrange the sequence so that there is no disorder in it. 

Consider the part of the sequence that starts at the beginning and ends right before a comparable pair of cubes appears. Then the cubes of this part belong to $M(Q)$, while the next cube does not. Continuing the process, we see that each lower subset of $Q$ can be constructed as follows. First we remove some cubes (call this set $R_1$) from $M(Q)=:M(Q_1)$, then we discard some set $R_2$ of cubes from $M(Q_1\setminus R_1)\setminus M(Q_1)=:M(Q_2)\setminus M(Q)$, and so on. In doing so, the number of ways to take away cubes in the first step is $\binom{|M(Q)|}{|R_1|}$, and in the $i$th step, for $i>1$, is $\binom{|M(Q_i)|-(|M(Q_{i-1})|-|R_{i-1}|)}{|R_i|}$. Denoting $k_i:=T(|Q_i|)-|M(Q_i)|+|R_i|,$ we have $\binom{|M(Q)|}{|R_1|}\leq \binom{T(|Q|)}{k_1}$ and for $i>1$,
\begin{align*}
\binom{|M(Q_i)|-(|M(Q_{i-1})|-|R_{i-1}|)}{|R_i|}&=\binom{|M(Q_i)|-(T(|Q_{i-1}|)-k_{i-1})}{|R_i|}\\
&\leq \binom{T(|Q_i|)-(T(|Q_{i-1}|)-k_{i-1})}{k_i}\leq \binom{k_{i-1}}{k_i}.
\end{align*}
Hence, the number of ways to construct a lower subset of $Q$ with a fixed sequence of $|R_i|$ is at most
\begin{align}\label{this}
&\binom{T(n)}{k_1}\binom{k_1}{k_2}\binom{k_2}{k_3}...\binom{k_{l-1}}{k_l}\leq \binom{T(n)}{k_1}2^{k_1+k_2+...+k_{l-1}}<\binom{T(n)}{k_1}2^{k},
\end{align}
where $l$ is the number of steps. If $T(|Q|)\leq k$, the right-hand side is bounded by $2^{2k}$. Otherwise, according to Stirling's formula,
\begin{align*}
\binom{T(n)}{k_1}2^{k}<\Big(\frac{eT(n)}{k_1}\Big)^{k_1}2^k\leq \Big(\frac{2eT(n)}{k}\Big)^{k}.
\end{align*}
Finally,
\begin{align}\label{333}
C(Q,k,d)\leq \sum_{|R_1|+...+|R_l|\leq k}\max\Big\{4,\frac{2eT(n)}{k}\Big\}^k<\Big(\max\Big\{8,\frac{4eT(n)}{k}\Big\}\Big)^k.
\end{align}
\end{proof}

\begin{corollary} If $n\geq d^{6d\log d}$, there holds
\begin{align}\label{cor1}
C(Q,k,d)<\Big(e^4\max\Big\{1,\frac{n^{1-\frac{1}{d}}}{k}\Big\}\Big)^k
\end{align}
and
\begin{align}\label{same}
C(Q,k,d)< 2^{2k+4n^{1-\frac{1}{d}}}.
\end{align}
\end{corollary}

\begin{proof}
Inequality \eqref{cor1} follows immediately from Lemmas \ref{lem} and \ref{cor} (note that $2\sinh \pi/\pi<e^2$). The second estimate \eqref{same} is valid due to Lemma \ref{lem} and relation \eqref{this} in the same fashion as \eqref{333}.
\end{proof}

Before proving Theorem \ref{pr}, let us establish a weaker result under less restrictive assumptions, which eventually will be helpful in the proof of Theorem \ref{pr}.

\begin{prop}\label{weaker}
For any $d\geq 2$ and $n\geq d^{12d\log d}$, we have
\begin{align*}
1<\frac{\log p_d(n)}{n^{1-\frac{1}{d}}}<d^2.
\end{align*}
\end{prop}

\begin{proof}
Let us first prove by induction on $d$ that 
\begin{align}\label{ind1}
\log p_{d}(n)< U_{d}n^{1-\frac{1}{d}},\qquad \text{for}\;n>d^{12d\log d},
\end{align}
where $U_2:=2\sqrt{2}$ and for $d\geq 3$,
\begin{align*}
U_d:=U_{d-1}d^{\frac{1}{d-1}}+1.
\end{align*}
The basis $d=2$ follows from the estimate $p_2(k)< e^{2\sqrt{2k}}$ (see at the end of \cite[Section 2]{HR2}). Assuming that \eqref{ind1} holds for $2,3,...,d-1$ we will prove it for $d\geq 3$. Take a lower set $Q,\;|Q|=n,$ set $k:=\lfloor n^{1/d}\rfloor$, and for a fixed $p,\;p=1,...,d,$ consider the following $``$slices$"$ of $Q$:
$$Q_i^p:=\{Q\cap \{q_p=i\}\}\setminus \bigcup\limits_{0< t<p,\;0\leq j< k} Q_j^t,\qquad i=0,...,k-1,$$ of cardinalities $n_0^p\geq n_1^p\geq ... \geq n_{k-1}^p$. 

For the sake of clarity, let us consider a three-dimensional example of such a slicing in Figure \ref{pic2}. 
\begin{figure}[H]
\center{\includegraphics[width=0.48\linewidth]{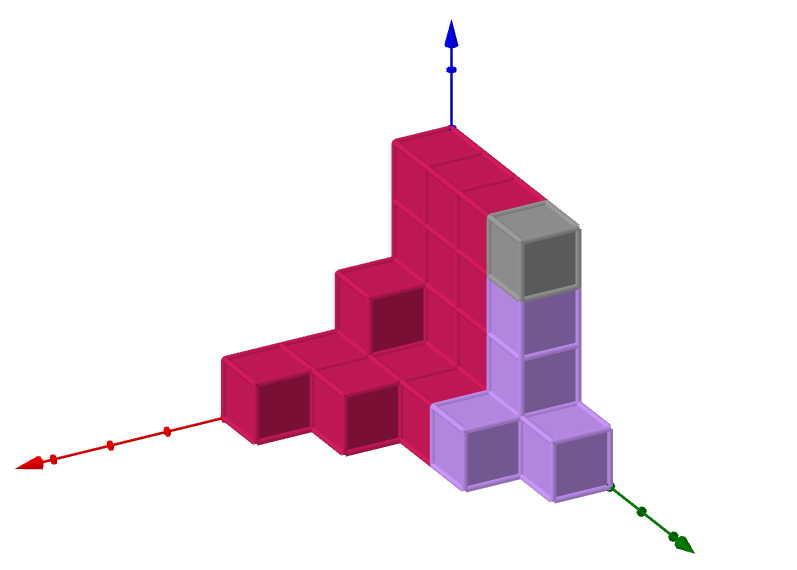}}
\caption{}
\label{pic2}
\end{figure}
Let us enumerate the axes in the following order: the green axis, the blue one, the red one. In our case we have $|Q|=25,\;k=\lfloor 25^{\frac{1}{3}}\rfloor=2$. Thus, $Q_0^1,Q_1^1,Q_2^1$ are the three magenta slices in Figure \ref{pic2} representing integer partitions of the numbers $|Q_0^1|=8,\;|Q_1^1|=6,$ and $|Q_2^1|=5$. The next slices $Q_0^2,Q_1^2,$ and $Q_2^2$ are the three violet layers with $|Q_0^2|=3,\;|Q_1^2|=1,$ and $|Q_2^2|=1$. Finally, the only nonempty slice along the third axis is $Q_0^3$ consisting of one gray cube.

Note that $Q=\bigcup_{0\leq i<k,\;1\leq p\leq d}Q_i^p$, since otherwise there would exist a cube $q\in Q$ with $q_i>n^{1/d}-1$ for all $i=1,...,d$, and, by the definition of lower sets, the cardinality of $Q$ would exceed $(n^{1/d})^d$, which is not true. So, 
$$\sum_{p=1}^d l_p=n,\quad l_p:=n_0^p+...+n_{k-1}^p.$$
In addition, any $Q_i^p,\;i>0,$ is a lower subset of $Q_{i-1}^p$, so once $Q_{i-1}^p$ is constructed, then if $n_i^p\geq (d-1)^{12(d-1)\log (d-1)}$, the number of possible $Q_i^p$ (with a fixed $n_i^p$) can be estimated either by \eqref{same} or by the induction assumption. As in \eqref{extra}, one can show that $(d-1)^{12(d-1)\log (d-1)}<d^{-6}n^{1-\frac{1}{d}}$. Thus, combining \eqref{cor1}, \eqref{same}, and the induction assumption we obtain the following bound for the logarithm of the number of slices of fixed cardinalities $n_i^p$ (for simplicity, we omit the upper indexes $p$ for $n_i^p$):
\begin{multline}\label{st1}
\max_{n_0+...+n_{k-1}=l_p}\Big\{U_{d-1}n_0^{1-\frac{1}{d-1}}+\sum_{i=0}^{k-2}\min\{2(n_i-n_{i+1})+4n_i^{1-\frac{1}{d-1}}, U_{d-1} n_{i+1}^{1-\frac{1}{d-1}}\}\Big\}
\\
+\sum_{i:\;n_i<d^{-6}n^{1-\frac{1}{d}}}(n_i-n_{i+1})(4+\log n_i)=:G(l_p,d)+\sum_{i:\;n_i<d^{-6}n^{1-\frac{1}{d}}}(n_i-n_{i+1})(4+\log n_i).
\end{multline}
Note that
\begin{align}\label{st3}
\sum_{i:\;n_i<d^{-6}n^{1-\frac{1}{d}}}(n_i-n_{i+1})(4+\log n_i)\leq\frac{n^{1-\frac{1}{d}}}{d^6}\big(4+\log (d-1)^{12(d-1)\log(d-1)}\big)<\frac{24n^{1-\frac{1}{d}}}{d^4}.
\end{align}
Then, taking into account \eqref{st1}, \eqref{st3}, and the inequality $p_2(n)\leq e^{2\sqrt{2n}}$, we derive
\begin{align*}
\log p_{d}(n)&\leq \log \Big|\Big\{\{n_i^p\},\;1\leq i\leq k,\;1\leq p\leq d:\sum_{i,p}n_i^p=n\Big\}\Big|+\sum_{p=1}^d G(l_p,d)+d\frac{24n^{1-\frac{1}{d}}}{d^4}\\
&\leq \log \binom{n+d-2}{d-1}+2d\sqrt{2n}+\frac{8n^{1-\frac{1}{d}}}{9}+\max_{\substack{l_1+...+l_d=n\\ n_0^p+...+n_{k-1}^p=l_p}}\sum_{p=1}^d \sum_{i=0}^{k-1}U_{d-1}n_{i}^{1-\frac{1}{d-1}}\\
&\leq d\log 2n+2d\sqrt{2n} +\frac{8n^{1-\frac{1}{d}}}{9}+U_{d-1} \max_{l_1+...+l_d=n}\sum_{p=1}^d n^{\frac{1}{d(d-1)}}l_p^{1-\frac{1}{d-1}}\\
&\leq 3d\sqrt{2n}+\frac{8n^{1-\frac{1}{d}}}{9}+U_{d-1}d^{\frac{1}{d-1}}n^{1-\frac{1}{d}}< U_d n^{1-\frac{1}{d}},
\end{align*}
completing the proof of \eqref{ind1}. 

Further, as $U_{d}<d^2$ for $d<8$, it suffices to prove the second part of the theorem by induction on $d\geq 8$. As above, the induction assumption holds for $Q_i^p$ with $n_i^p\geq (d-1)^{12(d-1)\log (d-1)}$. Likewise in \eqref{st1}, by \eqref{st3}, the logarithm of the number of slices of fixed cardinalities $n_i^p$ does not exceed
\begin{align*}
\tilde{G}(l_p,d)+\sum_{i:\;n_i<d^{-6}n^{1-\frac{1}{d}}}(n_i-n_{i+1})(4+\log n_i)\leq \tilde{G}(l_p,d)+\frac{24n^{1-\frac{1}{d}}}{d^4},
\end{align*}
with
\begin{align*}
\tilde{G}(l_p,d)&:= \max_{n_0+...+n_{k-1}=l_p}\Big\{(d-1)^2n_0^{1-\frac{1}{d-1}}+\sum_{i=0}^{k-2}\min\{2(n_i-n_{i+1})+4n_i^{1-\frac{1}{d-1}}, (d-1)^2n_{i+1}^{1-\frac{1}{d-1}}\}\Big\}\\
&\leq 4k\Big(\frac{l_p}{k}\Big)^{1-\frac{1}{d-1}}+\max_{n_0+...+n_{k-1}=l_p}\Big\{B_{d-1}n_0^{1-\frac{1}{d-1}}+\sum_{i=0}^{k-1}\min\{2(n_i-n_{i+1}), B_{d-1}n_{i+1}^{1-\frac{1}{d-1}}\}\Big\},
\end{align*}
where $B_{x}:=x^2-4,\;x\geq 8$. Noting that
\begin{align*}
B_{d-1}n_0^{1-\frac{1}{d-1}}&+\sum_{i=0}^{k-1}\min\{2(n_i-n_{i+1}), B_{d-1}n_{i+1}^{1-\frac{1}{d-1}}\}\\
&\leq B_{d-1}\sum_{i=0}^{\lfloor 2d^{-1}n^{1/d}\rfloor-1}n_i^{1-\frac{1}{d-1}}+2n_{\lfloor 2d^{-1}n^{1/d}\rfloor-1}\\
&\leq B_{d-1}2^{\frac{1}{d-1}}d^{-\frac{1}{d-1}}l_p^{1-\frac{1}{d}}+\frac{2l_p}{2d^{-1}n^{\frac{1}{d}}-1}< B_{d-1}2^{\frac{1}{d-1}}d^{-\frac{1}{d-1}}l_p^{1-\frac{1}{d}}+dl_pn^{-\frac{1}{d}}+2d^2l_pn^{-\frac{2}{d}}
\end{align*}
and that
\begin{align*}
4k\Big(\frac{l_p}{k}\Big)^{1-\frac{1}{d-1}}\leq 4n^{\frac{1}{d(d-1)}}l_p^{1-\frac{1}{d-1}},
\end{align*}
we have
\begin{align*}
\sum_{p=1}^d \tilde{G}(l_p,d)&\leq 4d^{\frac{1}{d-1}}n^{1-\frac{1}{d}}+B_{d-1}2^{\frac{1}{d-1}}n^{1-\frac{1}{d}}+dn^{1-\frac{1}{d}}+2d^2n^{1-\frac{2}{d}}.
\end{align*}

Hence, using again $p_2(n)\leq e^{2\sqrt{2n}}$, we obtain for $d\geq 8$
\begin{align*}
\log p_d(n)&\leq \log \Big|\Big\{\{n_i^p\},\;1\leq i\leq k,\;1\leq p\leq d:\sum_{i,p}n_i^p=n\Big\}\Big|+\sum_{p=1}^d \tilde{G}(l_p,d)+\frac{24n^{1-\frac{1}{d}}}{d^3}\\
&<3d\sqrt{2n}+4d^{\frac{1}{d-1}}n^{1-\frac{1}{d}}+B_{d-1}2^{\frac{1}{d-1}}n^{1-\frac{1}{d}}+dn^{1-\frac{1}{d}}+2d^2n^{1-\frac{2}{d}}+\frac{24n^{1-\frac{1}{d}}}{d^3}\\
&<3d\sqrt{2n}+2d^2n^{1-\frac{2}{d}}+\Big(d^{2}-\frac{d}{16}\Big)n^{1-\frac{1}{d}}<d^{2}n^{1-\frac{1}{d}},
\end{align*}
where in the last step we used that $n>\max\{(64d)^d,(96\sqrt{2})^6\}.$ This concludes the proof of the proposition.
\end{proof}

Now we are in a position to prove our main result.

\begin{proof}[Proof of Theorem \ref{pr}.] 
Our aim will be to prove the stronger inequality
\begin{align}\label{strong}
p_d(n)<1800K_{d}Y_dn^{1-\frac{1}{d}}=:CK_{d}Y_dn^{1-\frac{1}{d}},
\end{align}
where $Y_d:=\prod_{k=3}^{d-1}(1+2^{-k})^{1/k}$ and $K_d=\prod_{k=1}^{d-1}(1+k^{-2})$ as above. Inequality \eqref{strong} implies the needed one, since $K_d<\sinh \pi/\pi,\;\log Y_d<1/12$, and $e^{1/12}\sinh\pi/\pi<4.$

Let us prove \eqref{strong} by induction on $d$. The cases $d\leq 59$ follow from Proposition \ref{weaker}, so we prove the claim for $d\geq 60$ assuming that it holds for $2,3,...,d-1$. We divide the proof into several steps. 

{\bf Step 1. Partition into slices.} Take a lower set $Q,\;|Q|=n,$ and put
\begin{equation*}
k_p:=\begin{cases}
\lfloor 2^{-(d-1)}d^{-1}n^{1/d}\rfloor=:k,\qquad 0\leq p\leq d-1,\\
\lfloor 2^{d(d-1)} d^{d-1}n^{1/d}\rfloor, \qquad p=d.
\end{cases}
\end{equation*}
For a fixed $p,\;p=1,...,d,$ consider the slices $$Q_i^p:=\{Q\cap \{q_p=i\}\}\setminus \bigcup\limits_{0< t<p,\;0\leq j< k_t} Q_j^t,\qquad i=0,...,k_p-1$$ of cardinalities $n_0^p\geq n_1^p\geq ... \geq n_{k_p-1}^p$. Note that $\cup_{i,p}Q_i^p=Q$, since otherwise there exists a cube $q$ in $Q$ with $q_p\geq k_p$ for all $p$, so all the cubes with $p$th coordinate at most $k_p$ belong to $Q$ as well, which contradicts the fact that $\prod_p k_p>n$.

The idea is to split our lower set $Q$ into two subsets: the union of the chosen slices of the first $d-1$ directions $\bigcup_{1\leq p\leq d-1,\;0\leq i<k_p}Q_i^p$ and the complement subset $Q':=\bigcup_{0\leq i<k_d}Q_i^d$. Both $Q\setminus Q'$ and $Q'$ consist of slices that are lower sets themselves, and the number of lower subsets of the former can be well estimated using the induction assumption, as the number of slices is small enough. The remaining set $Q'$ has more slices, but we will see that the number of them is still well bounded, so that Lemma \ref{cor} can come into play $``$prohibiting$"$ the cardinalities of slices being close to each other. Denote $|Q'|=:l,\;|Q\setminus Q'|=:t$, so $l+t=n$. 
\medskip

{\bf Step 2. Dealing with small slices.} Observe that each of the chosen slices is a lower set and for each $p$ and $i<k_p-1$ the set $Q_{i+1}^p$ is a subset of $Q_i^p$. We can apply the induction assumption to $Q_i^p$ with $n_i^p\geq (30(d-1))^{2(d-1)^2}$. Noting that
\begin{align*}
(30(d-1))^{2(d-1)^2}<\Big(\frac{d}{d-1}\Big)^{2(d-1)^2}n^{\frac{(d-1)^2}{d^2}}\leq 2^{-2(d-1)}n^{\frac{(d-1)^2}{d^2}},
\end{align*}
for a fixed $p$ we have
\begin{align}\label{fixedp}
\sum_{n_i^p\leq 2^{-2(d-1)}n^{\frac{(d-1)^2}{d^2}}}(n_i^p-n_{i+1}^p)(4+\log n)<2^{-2d+3}n^{\frac{(d-1)^2}{d^2}}\log n.
\end{align}
\medskip

{\bf Step 3. Estimating the number of possible ${\bf Q\setminus Q'}$.} Using the induction assumption and the bound of the number of lower subsets given by \eqref{same} along with \eqref{fixedp}, we obtain that the logarithm of the number of possible $Q\setminus Q'$ with fixed $n_i^p$ is less than 
\begin{align}\label{first}
\sum_{p=1}^{d-1}&\sum_{i=0}^{k-1}CK_{d-1}Y_{d-1}(n_i^p)^{1-\frac{1}{d-1}}+(d-1)\cdot 2^{-2d+3}n^{\frac{(d-1)^2}{d^2}}\log n\nonumber\\
&< CK_{d-1}Y_{d-1}((d-1)k)^{\frac{1}{d-1}}t^{1-\frac{1}{d-1}}+n^{\frac{(d-1)^2}{d^2}}\log n\nonumber\\
&\leq 0.5CK_{d-1}Y_{d-1}n^{\frac{1}{d(d-1)}}t^{1-\frac{1}{d-1}}+n^{\frac{(d-1)^2}{d^2}}\log n.
\end{align}
\medskip

{\bf Step 4. Obtaining a general bound for the number of possible ${\bf Q'}$.} Now we estimate the number of possible $Q'$ with fixed $n_i:=n_i^d$. For the sake of simplicity let $m:=k_{d},\;\Delta_i:=n_i-n_{i+1},\;\Gamma_i=n_i^{1-1/(d-1)}$ (note that the last notation slightly differs from that of the proof of Lemma \ref{lem}). Combining the induction assumption with \eqref{cor1} and keeping in mind \eqref{fixedp}, we see that the logarithm of the number of lower sets with fixed $n_i$ cannot exceed the following sum (assuming $n_{m+1}=0$)
\begin{align}\label{rec}
CK_{d-1}Y_{d-1}&n_0^{1-\frac{1}{d-1}}+2^{-2d+3}n^{\frac{(d-1)^2}{d^2}}\log n+\sum_{i=0}^{m}\min\big\{(n_i-n_{i+1})\Big(4+\log^+\frac{n_i^{1-\frac{1}{d-1}}}{n_i-n_{i+1}}\Big),CK_{d-1}Y_{d-1}n_{i+1}^{1-\frac{1}{d-1}}\big\}\nonumber\\
&<4Cn_0^{1-\frac{1}{d-1}}+n^{\frac{(d-1)^2}{d^2}}\log n+\sum_{i=0}^{m}\min\big\{\Delta_i\Big(C+\log^+\frac{\Gamma_i}{\Delta_i}\Big), CK_{d-1}Y_{d-1}\Gamma_i\big\}.
\end{align}
The latter sum we can bound by
\begin{align}\label{rec1}
\sum_{i:\;\Delta_i\leq\Gamma_i}\Delta_i\log\frac{e^C\Gamma_i}{\Delta_i}+\sum_{i:\;\Delta_i> \Gamma_i}CK_{d-1}\Gamma_i=:\sum_{i=0}^{s-1} M_i+\sum_{i=s}^m M_i=:G(n_0,...,n_{s})+H(n_s,...,n_m),
\end{align}
where $s$ is the first index $i$ such that $n_i\leq l^{1-1/d}/4d^3.$ Note that in \eqref{rec} we can assume that $n_m\geq (30(d-1))^{2(d-1)^2}$, as the other $n_i$ are already taken into account.
\medskip

{\bf Step 5. Estimating ${\bf G(n_0,...,n_{s})}$.} Take a tuple $(n_0,...,n_{s})$ that delivers the maximum of the function $G$ over all the tuples in
 $$S'_l:=\Big\{(n_0,...,n_{s})\in\mathbb{Z}^+: n_0\geq ... \geq n_{s}\geq 0,\;n_{s-1}\geq \max\Big\{\frac{l^{1-\frac{1}{d}}}{4d^3},(30(d-1))^{2(d-1)^2}\Big\},\; n_0+...+n_{s}=l\Big\}.$$ 
Note that $\Delta_{s-1}> \Gamma_{s-1}$, since otherwise $n_s>0$ and we could decrease $n_s$ so that $\Delta_{s-1}> \Gamma_{s-1}$, which would increase the value of $G$.

Assume that for some $i,\;0\leq i\leq s-1,$ we have $$\Delta_i>\Gamma_i+2.$$ Then $\Delta_i> 3$ and if we substitute the pair $(n_i,n_{i+1})$ by $(n_i-1,n_{i+1}+1)$, the tuple will still be in $S'_n$ with $M_j,\;j\neq i-1,i,i+1,$ and $\Gamma_{i-1}$ unchanged. At the same time $\Delta_{i-1}$ increases, which means that $M_{i-1}$ does not decrease. Denote by $\Delta'_j, \Gamma'_j$ and $M'_j$ the corresponding values after the substitution. Consider the three cases.

Case 1. $\Delta_{i+1}> \Gamma_{i+1}$. Then $\Delta'_{i+1}> \Gamma'_{i+1}$ and 
\begin{align*}
(M'_i+M'_{i+1})-(M_i+M_{i+1})=CK_{d-1}((\Gamma'_i+\Gamma'_{i+1})-(\Gamma_i+\Gamma_{i+1}))>0.
\end{align*}

Case 2. $\Delta_{i+1}\leq \Gamma_{i+1},\; \Delta'_{i+1}\leq \Gamma'_{i+1}$. We have
\begin{multline*}
(M'_i+M'_{i+1})-(M_i+M_{i+1})\\
=(\Delta_{i+1}+1)\log \frac{e^C(n_{i+1}+1)^{1-\frac{1}{d-1}}}{\Delta_{i+1}+1}-\Delta_{i+1}\log \frac{e^Cn_{i+1}^{1-\frac{1}{d-1}}}{\Delta_{i+1}}+CK_{d-1}\big((n_i-1)^{1-\frac{1}{d-1}}-n_i^{1-\frac{1}{d-1}}\big)\\
\geq C-\Delta_{i+1}\log\frac{\Delta_{i+1}+1}{\Delta_{i+1}}-CK_{d-1}(0.5n_i)^{-\frac{1}{d-1}}>C-1-8Cn_i^{-\frac{1}{d-1}}>0,
\end{multline*}
since $n_i>16^{d-1}$.

Case 3. $\Delta_{i+1}\leq \Gamma_{i+1},\; \Delta'_{i+1}> \Gamma'_{i+1}$. Then
$$\Delta_{i+1}+1=\Delta'_{i+1}> \Gamma'_{i+1}>\Gamma_{i+1}>2,$$
so $\Delta_{i+1}\geq 0.5\Gamma_{i+1}$, and
\begin{multline*}
(M'_i+M'_{i+1})-(M_i+M_{i+1})=CK_{d-1}(n_{i+1}+1)^{1-\frac{1}{d-1}}-\Delta_{i+1}\log \frac{e^Cn_{i+1}^{1-\frac{1}{d-1}}}{\Delta_{i+1}}+CK_{d-1}\big((n_i-1)^{1-\frac{1}{d-1}}-n_i^{1-\frac{1}{d-1}}\big)\\
\geq CK_{d-1}(n_{i+1}+1)^{1-\frac{1}{d-1}}-\Delta_{i+1}\log 2e^C-CK_{d-1}(0.5n_i)^{-\frac{1}{d-1}}\\
\geq 0.25CK_{d-1}(1-4(0.5n_i)^{-\frac{1}{d-1}})>0,
\end{multline*}
since $n_i>4^d$.

Thus, in all the cases $G$ increases, and we come to a contradiction, which yields that 
\begin{align}\label{g'}
\Delta_{i}\leq \Gamma_{i}+2<1.25\Gamma_{i},\qquad 0\leq i\leq s-1,
\end{align}
as $n_i=\Gamma_i^{\frac{d-1}{d-2}}>8^{\frac{d-1}{d-2}}$.

Note that for $s\leq l^{1/d}$ we straightforwardly have an appropriate bound
\begin{align}\label{help}
G(n_0,...,n_s)\leq\max_{n_0+...n_{s-1}=l}\sum_{i=0}^{s-1}CK_{d-1}n_i^{1-\frac{1}{d-1}}= CK_{d-1}l^{1-\frac{1}{d}},
\end{align}
so from now on we assume 
\begin{align}\label{n_m'}
s-1\geq l^{\frac{1}{d}}\qquad\text{and}\qquad n_{s-1}\leq l^{1-\frac{1}{d}}.
\end{align} 
If $n_0> 2.5l^{1-1/d}$, then considering the sequence
\begin{align*}
a_0:=n_0\quad\text{ and}\quad a_{i}=a_{i-1}-1.25a_{i-1}^{1-\frac{1}{d-1}}\;\text{for}\;i\geq 1,
\end{align*}
similarly as in the proof of Lemma \ref{lem} (see Case a) we see that the sum of $a_i$ becomes equal to $l$ before $a_i$ reaches $l^{1-1/d}$, so the same holds for $n_i$, since $n_i$ decreases slower than $a_i$ (cf. \eqref{g'}). This contradicts \eqref{n_m'}. Thus, $$n_0\leq 2.5l^{1-\frac{1}{d}}.$$

Now, when we have the ratio $n_0/n_{s-1}$ bounded well enough by $10d^3$, we are going to show that $\Delta_{i}$ must be greater that $\Gamma_i$ for all $i=0,...,s-1$. Assume the contrary, so that for some $0<i<s$ there holds $\Delta_{i-1}\leq \Gamma_{i-1}$. Then we consider a new tuple $(n_0+n_i,n_1,...,n_{i-1},n_{i+1},...,n_s,0)$ instead of $(n_0,n_1,..,n_s)$ and estimate the difference between the values of $G$ at these points. First, let us estimate the difference $M'_0-M_0$.

Case a. $\Delta_0>\Gamma_0$. We have 
\begin{align*}
M'_0-M_0=CK_{d-1}((n_0+n_{i})^{1-\frac{1}{d-1}}-n_0^{1-\frac{1}{d-1}})\geq \frac{d-2}{d-1}CK_{d-1}\Gamma_i&\Big(\frac{n_{i}}{n_0+n_i}\Big)^{\frac{1}{d-1}}\\
&\geq CK_{d-1}\Gamma_i(10d^3)^{-\frac{1}{d-1}}=:W,
\end{align*}
as $n_{i}\geq l^{1-1/d}/4d^3$.

Case b. $\Delta_0\leq \Gamma_0,\;\Delta'_0\leq \Gamma'_0$. Then
\begin{align*}
M'_0-M_0&=(\Delta_0+n_{i})\log \frac{e^C(n_0+n_i)^{1-\frac{1}{d-1}}}{(\Delta_0+n_{i})}-\Delta_0\log \frac{e^Cn_0^{1-\frac{1}{d-1}}}{\Delta_0}\\
&\geq Cn_i-\Delta_0\log \frac{\Delta_0+n_i}{\Delta_0}\geq (C-1)n_i>W.
\end{align*}

Case c. $\Delta_0\leq \Gamma_0,\;\Delta'_0> \Gamma'_0$. In this case
\begin{align*}
M'_0-M_0&=CK_{d-1}(n_0+n_i)^{1-\frac{1}{d-1}}-\Delta_0\log \frac{e^Cn_0^{1-\frac{1}{d-1}}}{\Delta_0}\geq 2C(n_0+n_i)^{1-\frac{1}{d-1}}-C\Delta_0\frac{n_0^{1-\frac{1}{d-1}}}{\Delta_0}\geq W.
\end{align*}

Turn now to estimating the difference $M'_{i-1}-(M_{i-1}+M_i)$.

Case a'. $\Delta_i> \Gamma_i,\;\Delta_{i-1}+\Delta_{i}\leq \Gamma_{i-1}$. Then
\begin{align*}
M'_{i-1}-(M_{i-1}+M_i)&=(\Delta_{i-1}+\Delta_i)\log \frac{e^C\Gamma_{i-1}}{\Delta_{i-1}+\Delta_i}-\Delta_{i-1}\log \frac{e^C\Gamma_{i-1}}{\Delta_{i-1}}-CK_{d-1}\Gamma_i\\
&\geq C\Delta_i-\Delta_i-CK_{d-1}\Gamma_i\geq (C-1-CK_{d-1})\Gamma_i=:V.
\end{align*}

Case b'. $\Delta_i> \Gamma_i,\;\Delta_{i-1}+\Delta_{i}> \Gamma_{i-1}$. Note that
\begin{align*}
\Gamma_i=(\Gamma_{i-1}^{\frac{d-1}{d-2}}-\Delta_{i-1})^{\frac{d-2}{d-1}}>(\Gamma_{i-1}^{\frac{d-1}{d-2}}-\Gamma_{i-1})^{\frac{d-2}{d-1}}=\Gamma_{i-1}(1-\Gamma_{i-1}^{-\frac{1}{d-2}})^{\frac{d-2}{d-1}}\geq \frac{C}{C+1}\Gamma_{i-1},
\end{align*} 
since $n_{i-1}\geq (30(d-1))^{2(d-1)^2}> 60^{4(d-1)}>1801^{d-1}$. So,
\begin{align*}
M'_{i-1}-(M_{i-1}+M_i)&=CK_{d-1}\Gamma_{i-1}-\Delta_{i-1}\log \frac{e^C\Gamma_{i-1}}{\Delta_{i-1}}-CK_{d-1}\Gamma_i>-C\Gamma_{i-1}\geq (-1-C)\Gamma_{i}\geq V.
\end{align*}

Case c'. $\Delta_i\leq\Gamma_i$. We have
\begin{align*}
M'_{i-1}-(M_{i-1}+M_i)&=(\Delta_{i-1}+\Delta_i)\log \frac{e^C\Gamma_{i-1}}{\Delta_{i-1}+\Delta_i}-\Delta_{i-1}\log \frac{e^C\Gamma_{i-1}}{\Delta_{i-1}}-\Delta_{i}\log \frac{e^C\Gamma_i}{\Delta_{i}}\\
&\geq -\Delta_{i}\log \frac{e^C\Gamma_{i}}{\Delta_{i}}\geq -C\Gamma_{i}>V.
\end{align*}

Hence, in all the cases
\begin{align*}
G(n_0+n_i,...,n_m,0)-G(n_0,..,n_m)&\geq W+V=CK_{d-1}\Gamma_i(10d^3)^{-\frac{1}{d-1}}+(C-1-CK_{d-1})\Gamma_i\\
&>CK_{d-1}\Gamma_i\Big((10d^3)^{-\frac{1}{d-1}}-\frac{3C+1}{4C}\Big)>0
\end{align*}
for $d\geq 60$. This means that we come to a contradiction that ensures

\begin{align}\label{ha}
\Delta_{i}> \Gamma_i, \qquad 0\leq i\leq s-1.
\end{align}
With \eqref{ha} in hand, recalling also \eqref{help}, we can write
\begin{align}\label{g0}
G(n_0,...,n_s)&\leq\sum_{0\leq i<s:\;\Delta_i\leq\Gamma_i}\Delta_i\log\frac{e^C(\Delta_i+\Delta_{i+1})}{\Delta_i}+\sum_{0\leq i<s:\;\Delta_i> \Gamma_i}CK_{d-1}\Gamma_i\nonumber\\
&=\sum_{0\leq i<s:\;\Delta_i> \Gamma_i}CK_{d-1}\Gamma_i=\sum_{0\leq i<s:\;\Delta_i\leq\Gamma_i}C\Delta_i+\sum_{0\leq i<s:\;\Delta_i> \Gamma_i}CK_{d-1}\Gamma_i\nonumber\\
&<CK_{d-1}n^{1-\frac{1}{d}}\Big(1+\frac{2}{3(d-1)^2}+d^{-10(d-1)\log d}+l^{-\frac{1}{d}}\Big),
\end{align}
where the last inequality is due to \eqref{xi}.
\medskip

{\bf Step 6. Estimating ${\bf H(n_s,...,n_m)}$.} Let us split $H(n_s,...,n_m)$ (see \eqref{rec1}) into two sums
\begin{align*}
H(n_s,...,n_m)=\sum_{s\leq i:\;\Delta_i\leq 2^{-2d}d^{-8}\Gamma_i}M_i+\sum_{s\leq i:\;\Delta_i> 2^{-2d}d^{-8}\Gamma_i}M_i=:H_1+H_2,
\end{align*}
corresponding to, roughly speaking, big and small ratios $\Gamma_i/\Delta_i$. For $H_2$ and $d\geq 60$, we have the bound
\begin{align*}
H_2\leq n_s(4+\log 2^{2d}d^8)< 2dn_s\leq \frac{l^{1-\frac{1}{d}}}{2d^2}.
\end{align*}
Further, for $i$ satisfying $\Delta_i\leq 2^{-2d}d^{-8}\Gamma_i$, we obtain
\begin{align*}
\log\frac{e^4\Gamma_i}{\Delta_i}<e^2\sqrt{\frac{\Gamma_i}{\Delta_i}}\leq e^22^{-d}d^{-4}\frac{\Gamma_i}{\Delta_i},
\end{align*}
whence
\begin{align*}
H_1\leq \sum_{i=s}^{m}\frac{e^2}{2^dd^4}\Gamma_i\leq \frac{e^2}{2^dd^4} k_d^{\frac{1}{d-1}}l^{1-\frac{1}{d}}\leq \frac{e^2l^{1-\frac{1}{d}}}{d^3}.
\end{align*}
Thus,
\begin{align}\label{H}
H(n_s,...,n_m)=H_1+H_2\leq \frac{l^{1-\frac{1}{d}}}{2d^2}+\frac{e^2l^{1-\frac{1}{d}}}{d^3}<\frac{l^{1-\frac{1}{d}}}{d^2}.
\end{align}

{\bf Step 7. Combining all the estimates together.} Note that the number of different $n_i^p$ is less than $\binom{n+d-2}{d-1}(e^{2\sqrt{2n}})^d<e^{3d\sqrt{2n}}$. Therefore, recalling \eqref{first}, \eqref{rec}, \eqref{rec1}, \eqref{g0}, \eqref{H}, we infer
\begin{align*}
\log p_d(n)&\leq 0.5CK_{d-1}Y_{d-1}n^{\frac{1}{d(d-1)}}t^{1-\frac{1}{d-1}}+4Cn^{1-\frac{1}{d-1}}+2n^{\frac{(d-1)^2}{d^2}}\log n+3d\sqrt{2n}\\
&+CK_{d-1}Y_{d-1}l^{1-\frac{1}{d}}\Big(1+\frac{2}{3(d-1)^2}+d^{-10(d-1)\log d}+l^{-\frac{1}{d}}+\frac{1}{2Cd^2}\Big).
\end{align*}
Note that
\begin{align}\label{lii}
2n^{\frac{(d-1)^2}{d^2}}\log n<C n^{1-\frac{2}{d}+\frac{1}{d^2}}\log n\leq Cn^{1-\frac{2}{d}+\frac{1}{2d(d-1)}+\frac{1}{d^2}}<Cn^{1-\frac{1}{d-1}},
\end{align}
since $n^{\frac{1}{2d(d-1)}}\geq n^{\frac{\log\log n}{\log n}}=\log n$. Indeed, if the latter does not hold, then $d^2> \log n/2\log\log n$ and
\begin{align*}
d^{8d^2}> \Big(\frac{\log n}{2\log\log n}\Big)^{\frac{2\log n}{\log\log n}}\geq (\log n)^{\frac{\log n}{\log\log n}}=n,
\end{align*}
which contradicts the conditions of the theorem. So, using \eqref{lii} we come to
\begin{align}\label{fin0}
\log \;&p_d(n) <0.5CK_{d-1}Y_{d-1}n^{\frac{1}{d(d-1)}}t^{1-\frac{1}{d-1}}+5Cn^{1-\frac{1}{d-1}}+3d\sqrt{2n}+CK_{d-1}Y_{d-1}l^{1-\frac{1}{d}}\Big(1+\frac{3}{4(d-1)^2}+l^{-\frac{1}{d}}\Big)\nonumber\\
&<CK_{d-1}Y_{d-1}n^{\frac{1}{d(d-1)}}\Big(0.5 t^{1-\frac{1}{d-1}}+l^{1-\frac{1}{d-1}}\Big(1+\frac{3}{4(d-1)^2}+l^{-\frac{1}{d}}\Big)\Big)+5Cn^{1-\frac{1}{d-1}}+3d\sqrt{2n}.
\end{align}
Let us estimate the expression in brackets in the last line of \eqref{fin0}. If $l\leq n/2^d,$ then $$0.5t^{1-\frac{1}{d-1}}+l^{1-\frac{1}{d-1}}\Big(1+\frac{3}{4(d-1)^2}+l^{-\frac{1}{d}}\Big)<n^{1-\frac{1}{d-1}}.$$ 
Otherwise,
\begin{align*}
0.5t^{1-\frac{1}{d-1}}+l^{1-\frac{1}{d-1}}\Big(1+\frac{3}{4(d-1)^2}+l^{-\frac{1}{d}}\Big)<(0.5t^{1-\frac{1}{d-1}}+l^{1-\frac{1}{d-1}})\Big(1+\frac{3}{4(d-1)^2}+2n^{-\frac{1}{d}}\Big).
\end{align*}
Note that for any $0<a<b,\;\gamma\in (0,1),$ there holds
\begin{align*}
0.5a^{\gamma}+(b-a)^{\gamma}\leq b^{\gamma}(1+2^{-\frac{1}{1-\gamma}})^{1-\gamma},
\end{align*}
which in our case with $\gamma:=1-1/(d-1)$ gives
\begin{align*}
0.5t^{1-\frac{1}{d-1}}&+l^{1-\frac{1}{d-1}}\leq n^{1-\frac{1}{d-1}}(1+2^{-d+1})^{\frac{1}{d-1}}.
\end{align*}
Therefore in both cases we get
\begin{align}\label{fin1}
Y_{d-1}n^{\frac{1}{d(d-1)}}\Big(0.5 t^{1-\frac{1}{d-1}}+l^{1-\frac{1}{d-1}}\Big(1+\frac{3}{4(d-1)^2}+l^{-\frac{1}{d}}\Big)\Big)\leq Y_{d}n^{1-\frac{1}{d}}\Big(1+\frac{3}{4(d-1)^2}+2n^{-\frac{1}{d}}\Big).
\end{align}
Finally, \eqref{fin0} and \eqref{fin1} together give
\begin{align*}
\log p_d(n)&\leq CK_{d-1}Y_{d}n^{1-\frac{1}{d}}\Big(1+\frac{3}{4(d-1)^2}+2n^{-\frac{1}{d}}\Big)+5Cn^{1-\frac{1}{d-1}}+3d\sqrt{2n}< CK_{d}Y_{d}n^{1-\frac{1}{d}},
\end{align*}
where the latter inequality follows from
\begin{align*}
\max\{2n^{-\frac{1}{d}},2.5n^{-\frac{1}{d(d-1)}},dn^{\frac{1}{2}-1+\frac{1}{d}}\}=2.5n^{-\frac{1}{d(d-1)}}\leq \frac{1}{12(d-1)^2},
\end{align*}
as $n\geq (30d)^{2d^2}$. Thus, Theorem \ref{pr} is proved.
\end{proof}

\section{Lower sets in high dimensions}\label{sec2}

In cases of high dimensions the situation is rather different. In the first place, the trivial lower bound $p_d(n)\geq\binom{d+n-2}{d-1}$ becomes much more reasonable, since the configurations of lower sets in general become more sparse. We start with considering the case of a very large dimension $d$.

\begin{proof}[Proof of Theorem \ref{summing} (a).]
Put the first cube into the origin and for a fixed $j,\;0\leq j\leq n-1$, spread $j$ cubes along the axes. To complete a lower set, we have to add more $n-1-j$ cubes and we will do it stepwise. Note that any cube we put now is not situated along an axis, so it has at least two nonzero coordinates. This means that in any further step we must lean the current cube to at least two faces of some two cubes we put before. Since every pair of cubes can have at most one pair of their faces on which we can put a cube leaning, we come to the following estimate

\begin{align*}
p_d(n)\leq \sum_{j=0}^{n-1}\binom{d-1+j}{d-1}\prod_{k=j}^{n-2}\binom{k}{2}=:\sum_{j=0}^{n-1}A_j.
\end{align*} 
Noting that 
\begin{align*}
\frac{A_{j+1}}{A_j}=\frac{2(d+j)}{(j+1)j(j-1)}\geq\frac{2d}{n^3},
\end{align*}
we obtain
\begin{align*}
\frac{p_d(n)}{\binom{d+n-2}{d-1}}=
\frac{\sum_{j=0}^{n-1}A_j}{\binom{d+n-2}{d-1}}\leq \frac{1}{1-\frac{n^3}{2d}}\cdot\frac{A_{n-1}}{\binom{d+n-2}{d-1}}=\frac{1}{1-\frac{n^3}{2d}}.
\end{align*}
The estimate from below is given by 
\begin{align*}
p_d(n)\geq A_{n-1}=\binom{d+n-2}{d-1}.
\end{align*}
\end{proof}

The next result provides a more delicate estimate from above by dealing with a similar construction as in the proof of Theorem \ref{summing} (a). 

\begin{lemma}\label{main} 
There holds
\begin{align*}
p_d(n)\leq\sum_{m=2}^{n}\frac{e^m}{2^m}\sum_{t=1}^{m-1}(2\pi)^{-\frac{t+1}{2}}\sum_{\underset{s_i\geq 2,\;0\leq i<t}{s_0+...+s_t=m}}\frac{1}{\sqrt{s_0s_1...s_t}}(2d)^{s_0}s_0^{2s_1-s_0}s_1^{2s_2-s_1}...s_{t-1}^{2s_t-s_{t-1}}s_t^{-s_t}.
\end{align*}
\end{lemma}

\begin{proof}
Observe that every lower set can be constructed in the following way. First we put a cube into the origin. After that we choose some axes to put a cube along each of them, we call this zero step. Then, inductively, as we have completed the $(k-1)$th step, we have a lower set whose cubes have the sum of the coordinates less or equal to $k$. In the $k$th step we add some cubes to our set so that the following two conditions hold: any cube we put now has the sum of its coordinates equal to $k+1$ and the set we construct remains to be a lower set. 

Let us estimate the number of choices to put $s_k$ cubes in the $k$th step. Note that these $s_k$ cubes must lean only on $s_{k-1}$ cubes that we put in the previous step. When $k=1$, each pair of $s_0$ cubes from the previos step generates a place for a new cube and there are also $s_0$ possibilities to put a cube along an axis. So, the total number of possible places in this case is $s_0(s_0+1)/2<(s_0+1)^2/2$. Turn now to the cases of $k>1$.  Suppose that there are $l$ cubes among these $s_{k-1}$ ones that lie along some axes, that is, they have all the coordinates except one equal to zero. Then the only two ways to lean a new cube on any of these $l$ cubes are either to continue going along the corresponding axes or to lean it on one of these $l$ cubes and on one of the remaining $s_{k-1}-l$ ones. If a new cube does not lean on those $l$ cubes, then it has more than one nonzero coordinate, thus must lean on at least two cubes from the other $s_{k-1}-l$ ones from the previous step. As we have already noted, each pair of cubes generates at most one place for a new cube to lean on both of them. Summing up, the number of places to put cubes in the $k$th step is $1$ in the case $s_{k-1}=l=1$ and $$l+l(s_{k-1}-l)+\binom{s_{k-1}-l}{2}\leq \frac{s_{k-1}^2}{2},$$ otherwise. In the case $s_{k-1}=l=1$ all the remaining steps must have $s_i=1,\;i\geq k$. We come to the estimate
\begin{align*}
p_d(n)\leq\sum_{m=1}^{n-1}\sum_{t=1}^{m}\sum_{\underset{s_i\geq 2,\;1\leq i<t}{s_0+...+s_t=m}}\binom{d}{s_0}\binom{\frac{(s_0+1)^2}{2}}{s_1}\binom{\frac{s_1^2}{2}}{s_2}...\binom{\frac{s_{t-1}^2}{2}}{s_t}.
\end{align*}
Using Stirling's formula we see that
\begin{align*}
\binom{a}{b}\leq \frac{a^b}{b!}\leq \frac{1}{\sqrt{2\pi b}}\Big(\frac{ae}{b}\Big)^b
\end{align*}
for any $a\geq b\geq 1$, so we finally obtain
\begin{align*}
p_d(n)&\leq \sum_{m=1}^{n-1}\sum_{t=1}^{m}\sum_{\underset{s_i\geq 2,\;1\leq i<t}{s_0+...+s_t=m}}e^{s_0}\Big(\frac{e}{2}\Big)^{m-s_0}(2\pi)^{-\frac{t+1}{2}}\frac{1}{\sqrt{s_0s_1s_2...s_{t}}}d^{s_0}s_0^{-s_0}(s_0+1)^{2s_1}s_1^{-s_1}s_1^{2s_2}s_2^{-s_2}...s_{t-1}^{2s_t}s_t^{-s_t}\\
&=\sum_{m=2}^{n}\frac{e^m}{2^{m}}\sum_{t=1}^{m-1}(2\pi)^{-\frac{t+1}{2}}\sum_{\underset{s_i\geq 2,\;0\leq i<t}{s_0+...+s_t=m}}\frac{1}{\sqrt{s_0s_1...s_t}}(2d)^{s_0}s_0^{2s_1-s_0}s_1^{2s_2-s_1}...s_{t-1}^{2s_t-s_{t-1}}s_t^{-s_t}.
\end{align*}
\end{proof}

Lemma \ref{main} will be our main tool for further upper estimates of $p_d(n)$. The first one is of interest when $d=o(n^2)$ as $n\to \infty$.

\begin{prop}\label{up}
For $d\leq n^2/4$, there holds
\begin{align*}
p_d(n)< 4e^{cn}n^{n+2\sqrt{d}}\max\{2^{-n}, (2n)^{-\sqrt{d}}\}\quad \text{with}\;\;c=\frac{3}{2e}+1,
\end{align*}
which in case $dn^{-2}\to 0$ as $n\to\infty$, yields
\begin{align*}
p_d(n)\leq n^{n+o(n)}.
\end{align*}
\end{prop}

\begin{proof}
Consider a tuple $(s_0,s_1,...,s_t)\in \mathbb{N}^{t+1}$ such that $s_i\geq 2$ for $0<i<t$ and $s_0+...+s_t=m$. Note that 
\begin{align}\label{fir}
(s_0s_1...s_t)^{\frac{3}{2}}\leq \Big(\frac{m}{t}\Big)^{\frac{3t}{2}}\leq e^{\frac{3m}{2e}}.
\end{align}
Suppose that there is no $s_i,\;i\geq 1,$ such that $s_i\geq \sqrt{d}$. Then, using \eqref{fir}, we have
\begin{align}\label{lesssq}
F(d,s_0,...,s_t):&=(s_0s_1...s_t)^{\frac{3}{2}}(2d)^{s_0}s_0^{2s_1-s_0}s_1^{2s_2-s_1}...s_{t-1}^{2s_t-s_{t-1}}s_t^{-s_t}\leq e^{\frac{3m}{2e}}(2d)^{s_0}s_0^{2s_1-s_0}s_1^{s_2}s_2^{s_3}...s_{t-1}^{s_t}s_t^{-s_1}\nonumber\\
&\leq e^{\frac{3m}{2e}}(2d)^{s_0}s_0^{2\sqrt{d}-s_0}d^{\frac{m-s_0}{2}}=e^{\frac{3m}{2e}}d^{\frac{m+s_0}{2}}2^{s_0}s_0^{2\sqrt{d}-s_0}<e^{\frac{3m}{2e}}d^{\frac{m}{2}}m^{2\sqrt{d}}2^{s_0}d^{\frac{s_0}{2}}s_0^{-s_0}.
\end{align}
Here we used the inequality $s_1^{s_1}s_2^{s_2}...s_t^{s_t}\geq s_1^{s_2}s_2^{s_3}...s_t^{s_1}$, which is true since for any $k\in\mathbb{N}$, any positive integers $a_1,...,a_k$, and any permutation $\sigma$ from the symmetric group $\mathfrak{S}_k$, there holds 
\begin{align}\label{perm}
\prod_i a_i^{a_{\sigma(i)}}\leq\prod_i a_i^{a_i}.
\end{align} 
The maximum of the right-hand side of \eqref{lesssq} is attained at $s_0=2\sqrt{d}/e$, so in this case

\begin{align}\label{case1}
F(d,s_0,...,s_t)e^{-\frac{3m}{2e}}\leq d^{\frac{m}{2}}(e^{\frac{1}{e}}m)^{2\sqrt{d}}\leq n^{n}2^{-n}e^{\frac{n}{2}}n^{2\sqrt{d}}<n^{n+2\sqrt{d}}.
\end{align}
If there exists $s_i\geq\sqrt{d},\;i\geq 1,$ then choosing the maximal such index $i$ and using twice inequality \eqref{perm} along with \eqref{fir} we obtain
\begin{align}\label{case2}
F(d,s_0,...,s_t)e^{-\frac{3m}{2e}}&\leq 2^{s_0}(d^{s_0}s_0^{2s_1-s_0}...s_{i-1}^{2s_i-s_{i-1}}s_i^{-s_i})m^{2s_{i+1}}(s_{i+1}^{2s_{i+2}-s_{i+1}}...s_{t-1}^{2s_t-s_{t-1}}s_t^{-s_t})\nonumber\\
&\leq 2^{s_0}(s_0^{2s_1-s_0}...s_{i-1}^{2s_i-s_{i-1}}s_i^{2s_0-s_i})m^{2s_{i+1}}(s_{i+1}^{s_{i+2}}...s_{t-1}^{s_t}s_t^{-s_{i+1}})\nonumber\\
&\leq 2^{s_0}(s_0^{s_1}...s_{i-1}^{s_i}s_i^{s_0})m^{2s_{i+1}}m^{s_{i+2}+...+s_t}\nonumber\\
&\leq 2^{s_0}m^{s_0+...+s_i+2s_{i+1}+s_{i+2}+...+s_t}\nonumber\\
&<2^{m-\sqrt{d}} m^{m+\sqrt{d}}\leq 2^{n-\sqrt{d}} n^{n+\sqrt{d}}.
\end{align}
According to Lemma \ref{main} it remains only to estimate the sum
\begin{align*}
\sum_{m=2}^{n}\frac{e^m}{2^m} \sum_{t=1}^{m-1}(2\pi)^{-\frac{t+1}{2}}\sum_{\underset{s_i\geq 2,\;0\leq i<t}{s_0+...+s_t=m}}\frac{1}{(s_0s_1...s_t)^2}.
\end{align*}
Note that the right sum is equal to the coefficient at $x^m$ of the polynomial
\begin{align*}
P(x):=\Big(\sum_{j=1}^m \frac{x^j}{j^2}\Big) \Big(\sum_{j=2}^m \frac{x^j}{j^2}\Big)^{t}.
\end{align*}
We have $P(1)<(\pi^2/6)^{t+1}$, therefore,
\begin{align}\label{bas}
\sum_{m=2}^{n}\frac{e^m}{2^m} \sum_{t=1}^{m-1}(2\pi)^{-\frac{t+1}{2}}\sum_{\underset{s_i\geq 2,\;0\leq i<t}{s_0+...+s_t=m}}\frac{1}{(s_0s_1...s_t)^2}&<\sum_{m=2}^{n}\frac{e^m}{2^m} \sum_{t=1}^{m-1}\Big(\frac{\pi^2}{6\sqrt{2\pi}}\Big)^{t+1}\nonumber\\
&<\frac{4}{9}\frac{1}{1-\frac{2}{3}}\sum_{m=2}^{n}\frac{e^m}{2^m}=\frac{4}{3}\sum_{m=2}^{n}\frac{e^m}{2^m}.
\end{align}
Thus, combining \eqref{bas} with \eqref{case1} and \eqref{case2}, we finally derive
\begin{align*}
p_d(n+1)\leq n^{n+2\sqrt{d}}\max\bigg\{1, \frac{2^n}{(2n)^{\sqrt{d}}}\bigg\}\frac{4}{3}\sum_{m=1}^{n}\Big(\frac{e^{\frac{3}{2e}+1}}{2}\Big)^m<4e^{\frac{3n}{2e}+n}n^{n+2\sqrt{d}}\max\{2^{-n}, (2n)^{-\sqrt{d}}\},
\end{align*}
which concludes the proof.
\end{proof}

The complementary lower bound for the case $d=o(n^2)$ will be as follows.

\begin{prop}\label{low}
If $d\geq n/\psi(n)$ for some positive function $\psi(n)\geq 1$ such that $\psi(n)\to\infty$ and $\log \psi(n)/\log n\to 0$ as $n\to\infty$, then there holds
\begin{align*}
p_d(n+1)\geq n^{\big(n+\frac{n\log d}{\psi(n)\log n}\big)(1+o(1))},
\end{align*}
and consequently, if $\frac{\log d}{\psi(n)\log n}\to 0$ as $n\to\infty$, we have
\begin{align*}
p_d(n)\geq n^{n+o(n)}.
\end{align*}
\end{prop}

\begin{proof}
Let us count only the lower sets whose cubes have at most two nonzero coordinates and these coordinates are at most $1$. So first we put a cube into the origin, then we fix some number $i,\;i=1,...,n,$ and choose $i$ axes to put cubes along them. The remaining cubes will lie in the obvious way on some of the two-dimensional hyperplanes generated by the chosen axes. So the number of such lower sets is
\begin{align*}
\sum_{i=1}^{\min\{d,n\}}\binom{d}{i}\binom{\frac{i(i-1)}{2}}{n-i}=:\sum_{i=1}^{\min\{d,n\}}B_i\geq B_{\lfloor n/\psi(n)\rfloor}.
\end{align*}
To estimate the latter value, we note that by Stirling's formula
\begin{align*}
\sqrt{2\pi a}\Big(\frac{a}{e}\Big)^a\leq a!\leq e^{\frac{1}{12}}\sqrt{2\pi a}\Big(\frac{a}{e}\Big)^a
\end{align*}
for all $a\geq 1$, so this implies the inequality
\begin{align*}
\binom{a}{b}\geq \frac{\sqrt{a}a^a}{e^{\frac{1}{6}}\sqrt{2\pi(a-b)b}(a-b)^{a-b}b^b}>\frac{1}{2\sqrt{a}}\Big(\frac{a}{b}\Big)^b
\end{align*}
for all $a> b\geq 1$. Now, assuming that $\psi(n)\leq n/6$, we can estimate $B_{\lfloor n/\psi(n)\rfloor}$ as follows (writing just $\psi$ in place of $\psi(n)$ for the sake of simplicity)
\begin{align}\label{low_a}
&B_{\lfloor n/\psi\rfloor}=\binom{d}{\lfloor \frac{n}{\psi}\rfloor}\binom{\frac{1}{2}\lfloor \frac{n}{\psi}\rfloor ^2-\frac{1}{2}\lfloor \frac{n}{\psi}\rfloor }{n-\lfloor \frac{n}{\psi}\rfloor}\nonumber\\
&> \frac{1}{2\sqrt{d}}\Big(\frac{d\psi}{n}\Big)^{\frac{n}{\psi}-1}\frac{\psi}{\sqrt{2}n}\Big(\frac{n}{4\psi^2}\Big)^{n-\frac{n}{\psi}}>d^{\frac{n}{\psi}-\frac{3}{2}}n^{n-2\frac{n}{\psi}}\psi^{-2n+3\frac{n}{\psi(n)}}2^{-2n+\frac{2n}{\psi}-\frac{3}{2}}\nonumber\\
&=\exp\Big(n\log n-\frac{2n\log n}{\psi}+\frac{n\log d}{\psi}-\frac{3\log d}{2}-2n\log \psi+\frac{3n\log\psi}{\psi}-2n\log 2+\frac{2n\log 2}{\psi}-\frac{3\log 2}{2}\Big)\\
&=n^{\big(n+\frac{n\log d}{\psi\log n}\big)(1+o(1))},\nonumber
\end{align}
which in case $\log d=o(\psi(n)\log n)$ yields
\begin{align*}
B_{\lfloor n/\psi(n)\rfloor}\geq n^{n+o(n)}.
\end{align*}
The claim follows immediately from these estimates.
\end{proof}

\begin{remark} For $d\geq n/\log n$ and $\psi(n):=\log n$, inequality \eqref{low_a} gives
\begin{align*}
p_d(n+1)>n^{n-\frac{6n\log\log n}{\log^2 n}}.
\end{align*}
\end{remark}

\begin{proof}[Proof of Theorem \ref{summing} (d).]
The relation follows straightforwardly from the corresponding parts of Propositions \ref{up} and \ref{low}.
\end{proof}

Now we give a more general estimate, which will imply the sharp exponential order of $p_d(n)$ in case of $n^2=O(d)$.

\begin{prop}\label{4} If $d\geq \xi n^2$ for some $\xi=\xi(n)\geq 2n^{-1}$, then
\begin{align*}
p_d(n)< 3a^{2n}e^{\frac{125n}{\xi}}n^{2}e^{n}\frac{d^n}{n^n}\quad\text{with}\;\;a=\max\{2e^{3.5}\xi^{-1},1\}.
\end{align*}
In particular, if $\xi(n)\to\infty$ as $n\to\infty$, then 
\begin{align*}
p_d(n)= e^{n}\frac{d^n}{n^n}e^{o(n)}.
\end{align*}
\end{prop}

\begin{proof}
For a tuple $(s_0,...,s_t)\in\mathbb{N}^{t+1}$ with $s_0+...+s_t=m$, as before let
\begin{align*}
F(d,s_0,...,s_t):=(s_0s_1...s_t)^{\frac{3}{2}}(2d)^{s_0}s_0^{2s_1-s_0}s_1^{2s_2-s_1}...s_{t-1}^{2s_t-s_{t-1}}s_t^{-s_t}.
\end{align*}
We will prove by induction on $t$ that 
\begin{align}\label{ine}
F(d,s_0,...,s_t)\leq a^{2m-s_0}\exp\bigg(-\int\limits_{\min\{\frac{125m}{\xi},s_1\}}^{\frac{125m}{\xi}}\log\frac{\xi x}{125m}\;dx\bigg)m^{2}\frac{(2d)^m}{m^m}.
\end{align}
Note that for any $\alpha\leq\beta$ and $\gamma$,
\begin{align}\label{easy}
W(\a,\beta,\g):=\int\limits_{\g\a}^{\g\beta}\log \frac{x}{\g}\;dx=\g(x\log x-x)\Big|_{\a}^{\beta}.
\end{align}
The case $t=0$ is clear (we assume $s_1=0$). In the case $t=1$ we have
\begin{align*}
(\log F(d,s_0,m-s_0))'_{s_0}&=\Big( \log\Big((2d)^{s_0}s_0^{2m-3s_0+1.5}(m-s_0)^{-m+s_0+1.5}\Big)\Big)'_{s_0}=\log\frac{2d(m-s_0)}{s_0^3}\\
&+\frac{2(m-s_0)+1.5}{s_0}-\frac{1.5}{m-s_0}>\log\frac{2\xi(m-s_0)}{e^{1.5}m}>\log\frac{\xi(m-s_0)}{3m}.
\end{align*}
If $s_1\geq 6m/\xi,$ then applying the inequality above and \eqref{easy}, we see that
\begin{align*}
F(d,s_0,s_1)\leq \exp\big(-W(0,2,3m\xi^{-1})\big)F(d,m,0)\leq F(d,m,0)=\frac{(2d)^m}{m^m},
\end{align*}
which is less then the right-hand side of \eqref{ine}. Otherwise, $s_1< 6m/\xi<125m/2\xi$ and
\begin{align*}
F(d,s_0,s_1)\leq \exp\big(-W(0,1,3m\xi^{-1})\big)F(d,m,0)=e^{\frac{3m}{\xi}}\frac{(2d)^m}{m^m},
\end{align*}
while at the right-hand side of \eqref{ine} we get at least
\begin{align*}
\exp\big(-W(0.5,1,125m\xi^{-1})\big)\frac{(2d)^m}{m^m}=e^{\frac{125m}{2\xi}}\frac{(2d)^m}{m^m},
\end{align*}
which completes the proof of \eqref{ine} for $t=1$. 

Assume now that $t>1$ and \eqref{ine} is proved for all $m$ and for $1,...,t-1$. Let us prove it for $t$. The idea of the proof is to perform a chain of small perturbations of a tuple $(s_0,s_1,...,s_t)\in \mathbb{N}^{t+1},\;s_0+...+s_t=m,$ in such a way that at the end of this process the function $F(d,s_0,s_1,...,s_t)$ accumulates an insignificant factor and either we arrive at a shorter tuple and use the induction assumption or we come to a convenient case $s_1\geq 2s_2\geq 2^{i-1}s_i,\;s_{i+1}=...=s_t=1,$ for some $1\leq i\leq t$.

Consider a tuple $(s_0,s_1,...,s_t)\in \mathbb{N}^{t+1}$ such that $s_0+...+s_t=m$ and suppose that $s_t>s_{t-1}/2$. Fix $s_1,...,s_{t-1}$ and $s_0+s_t=:y$ and see what occures if we increase $s_0=:x$. We have
\begin{align*}
\Big(&\log F(d,x,s_1,...,s_{t-1},y-x)\Big)'_x \\
&=\big(x\log 2d+(2s_1-x+1.5)\log x+(2y-2x-s_{t-1}+1.5)\log s_{t-1}+(x-y+1.5)\log (y-x)\big)'_x\\
&>\log\frac{2d(y-x)}{xs_{t-1}^2}-\frac{1.5}{y-x}\geq\log\frac{2ds_t}{s_0s_{t-1}^2}-1.5>\log\frac{d}{s_0s_{t-1}}-1.5>\log\frac{d}{m^2}-1.5\geq\log \xi-1.5.
\end{align*}
This means that we can increase $s_0$ by $1$ and decrease $s_t$ keeping their sum constant until either $s_t\leq s_{t-1}/2$ or $s_t=1$ so that in every step of this process the value of $F(d,s_0,...,s_t)$ changes by at least $\exp(\log\xi-1.5)= \xi e^{-1.5}$. 

Suppose now that for some $i,\;1< i\leq t-1$ and $j\geq i$, we have 
\begin{align}\label{cond}
s_l\geq 2s_{l+1}\; \text{for}\;\;i\leq l<j\quad\text{and} \quad s_j=...=s_t=1.
\end{align}
Fix $s_1,...,s_{i-1},s_{i+1},...,s_t$ and $s_0+s_i=:y$ and observe what occures if we start increasing $s_0=:x$. We see that
\begin{align*}
&\Big(\log F(d,x,s_1,...,s_{i-1},y-x,s_{i+1},...,s_t)\Big)'_x \\
&\quad =\big(x\log 2d+(2s_1-x+1.5)\log x+(2y-2x)\log s_{i-1}+(2s_{i+1}+x-y+1.5)\log (y-x)\big)'_x\\
&\quad\quad >\log\frac{2d(y-x)}{xs_{i-1}^2}-\frac{2s_{i+1}+1.5}{y-x}=\log\frac{2ds_i}{s_0s_{i-1}^2}-\frac{2s_{i+1}+1.5}{s_i}\geq \log \xi-3.5,
\end{align*}
while $s_i\geq s_{i+1}$ (which is true under \eqref{cond}) and $s_i\geq s_{i-1}/2$. So we can decrease $s_i$ and increase $s_0$ with their sum constant, changing $F(d,s_0,...,s_t)$ by at least $\xi e^{-3.5}$ in each step, until one of the following situations happens.

Case 1. $s_i=s_{i+1}=...=s_t=1$. Then we accumulated at most the extra factor $(e^{3.5}\xi^{-1})^{\Delta s_0}$, where by $\Delta s_0$ we denote the number of steps we made increasing $s_0$, which is exactly the difference between the value of $s_0$ in the end and in the beginning of the process. So, we come to \eqref{cond} with $i-1$ in place of $i$ and proceed inductively with $i-1$ instead of $i$.

Case 2. $s_{i-1}\geq 2s_i\geq 4s_{i+1}$. Then we come again to \eqref{cond} with $i-1$ in place of $i$ with the same accumulated factor.

Case 3. $s_{i-1}< 2s_i= 4s_{i+1}$. Then $s_i=2s_{i+1}=:2x$ and we merge $s_i=2x$ and $s_{i+1}=x$ into one single variable equal to $2x$. Let us compare the new value of $F$ with the original one: 
\begin{align*}
\frac{F(s_0,...,s_{i-1},2x,s_{i+2},...,s_t)}{F(s_0,...,s_t)}=\frac{s_{i-1}^{4x-s_{i-1}+1.5}(2x)^{2s_{i+2}-2x+1.5}}{s_{i-1}^{4x-s_{i-1}+1.5}(2x)^{2x-2x+1.5}x^{2s_{i+2}-x+1.5}}\geq (4x)^{-x-1.5}.
\end{align*}
At the same time the sum of all the variables $s_0,s_1,...,s_{i-1},2x,s_{i+2},...,s_t$ becomes $m-x$ instead of $m$. So, by the induction assumption we have
\begin{align*}
&F(d,s_0,...,s_t) \leq \max_{x}\;(4x)^{x+1.5}a^{2m-2x-s_0}\exp\bigg(-\int\limits_{\min\{\frac{125(m-x)}{\xi},s_1\}}^{\frac{125(m-x)}{\xi}}\log\frac{\xi y}{125(m-x)}\;dy\bigg)\\
&\times(m-x)^{2}\frac{(2d)^{m-x}}{(m-x)^{m-x}}\leq \max_{x}\;(4x)^{x+1.5}a^{2m-2x-s_0}\exp\bigg(-\int\limits_{\min\{\frac{125m}{\xi},s_1\}}^{\frac{125m}{\xi}}\log\frac{\xi y}{125m}\;dy\bigg)m^{2}\frac{(2d)^{m-x}}{(m-x)^{m-x}}.
\end{align*}
As
\begin{align*}
\Big((x+1.5)\log 4x+(m-x)(\log 2d-\log(m-x))\Big)'_x&=\log 4x+1-\log 2d+\frac{1.5}{x} +\log (m-x)+1\\
&=\log \frac{4e^{3.5}x(m-x)}{2d}<\log \frac{2e^{3.5}}{\xi},
\end{align*}
we finally have the desired inequality \eqref{ine}, since
\begin{align*}
F(d,s_0,...,s_t)&\leq \max_{x}\Big(\max\Big\{\frac{2e^{3.5}}{\xi},1\Big\}\Big)^x a^{2m-2x-s_0}\exp\bigg(-\int\limits_{\min\{\frac{125m}{\xi},s_1\}}^{\frac{125m}{\xi}}\log\frac{\xi y}{125m}\;dy\bigg)m^{2}\frac{(2d)^{m}}{m^{m}}\\
&=a^{2m-s_0}\exp\bigg(-\int\limits_{\min\{\frac{125m}{\xi},s_1\}}^{\frac{125m}{\xi}}\log\frac{\xi y}{125m}\;dy\bigg)m^{2}\frac{(2d)^{m}}{m^{m}}.
\end{align*}

This way we either merged two variables in some step and obtained the needed inequality using the induction assumption or came to the situation $s_1\geq 2s_2\geq 2^{i-1}s_i,\;s_{i+1}=...=s_t=1,$ for some $1\leq i\leq t$. In the latter occasion, considering $s_0+s_1=:y$ to be constant and changing $s_0=:x$ we see that
\begin{align*}
\Big(\log F(d,x,y-x,s_2,...,s_t)\Big)'_x &=(x\log 2d+(2y-3x+1.5)\log x+(2s_2-y+x+1.5)\log (y-x))'_x\\
&\geq\log\frac{2d(y-x)}{x^3}-2-\frac{2s_2+1.5}{y-x}\geq\log\frac{2ds_1}{s_0^3}-\frac{2s_2}{s_1}-3.5.
\end{align*}
Thus, while $s_1\geq s_2$, there holds
\begin{align*}
\Big(\log F(d,x,y-x,s_2,...,s_t)\Big)'_x \geq\log\frac{2ds_1}{e^{5.5}s_0^3}>\log\frac{\xi s_1}{125m}.
\end{align*}
So, we can decrease $s_1$ and increase $s_0$ with $s_0+s_1$ constant so that $F(d,s_0,...,s_t)$ in this process increases at least by $\exp(\int_{s_1'}^{s_1^*}\log(\xi x/125m)\;dx)$, where $s_1^*$ stands for the value of $s_1$ that we started from and $s_1'$, for the value where we stopped, until one of the following situations happens.

Case a. $s_1=s_2=...=s_t=1$. Then 
\begin{align*}
F(d,s_0,...,s_t)&\leq \exp\bigg(-\int\limits_{s_1}^{s_1^*}\log\frac{\xi x}{125m}\;dx\bigg)\max_{s_0}\frac{s_0^2(2d)^{s_0}}{s_0^{s_0-2}}\leq \exp\bigg(-\int\limits_{\min\{\frac{125m}{\xi},s_1\}}^{\frac{125m}{\xi}}\log\frac{\xi x}{125m}\;dx\bigg)\frac{(2d)^m}{m^{m-2}},
\end{align*}
and the needed inequality is proved.

Case b. $s_1= 2s_2$. Then we merge $s_1$ and $s_2$ into $s_1$ as above and use the induction assumption. The only difference is that we have to take into account the factor $\exp(\int_{s_1'}^{s_1^*}\log(\xi x/125m)\;dx)$ that we accumulated while making $s_1$ decrease.

Thus, in all cases we obtained \eqref{ine} for all $m$.

So we have
\begin{align*}
F(d,s_0,...,s_t)\leq a^{2m}\exp\bigg(-\int\limits_{0}^{\frac{125m}{\xi}}\log\frac{\xi x}{125m}\;dx\bigg)m^{2}\frac{(2d)^m}{m^m},
\end{align*}
which in light of equality \eqref{easy} implies
\begin{align*}
F(d,s_0,...,s_t)\leq a^{2m}e^{\frac{125m}{\xi}}m^{2}\frac{(2d)^m}{m^m}.
\end{align*}
Finally, taking into account Lemma \ref{main} and estimate \eqref{bas}, we obtain
\begin{align*}
p_d(n)<a^{2n}e^{\frac{125n}{\xi}}n^{2}\frac{d^n}{n^n}\frac{4}{3}\sum_{m=2}^{n}e^m<3a^{2n}e^{\frac{125n}{\xi}}n^{2}\frac{d^n}{n^n}e^{n},
\end{align*}
which concludes the proof.
\end{proof}

\begin{proof}[Proof of Theorem \ref{summing} (b), (c).]
The claim follows from Proposition \ref{4} and the simple estimate $p_d(n)\geq \binom{d+n-2}{d-1}$.
\end{proof}

The following corollary combines the results above and applies them to the power--logarithmic scale of $d$ in terms of $n$.

\begin{corollary} If $cn^{\a}\log^{\g}n\leq d\leq Cn^{\a}\log^{\g}n$ for some $\a\geq 1,\;\g\in\mathbb{R}$, and positive constants $c$ and $C$, then
\begin{align*}
p_d(n)=\begin{cases} \binom{d+n-2}{d-1}\theta(d,n), &\quad\text{if }\; \a>3,\; \text{or}\;\; \a= 3,\;\g> 0,\\
e^{n}n^{(\a-1)n}\log^{\g n}n \;e^{O(n^{3-\a}\log^{-\g}n+\log n)}, &\quad\text{if }\; 2\leq \a\leq 3,\\
n^ne^{O(n\log\log n)}, &\quad\text{if }\;1\leq\a<2.
\end{cases}
\end{align*}
Here the function $\theta (d,n)\geq 1$ is bounded from above by a constant that depends only on $\a$ and $\g$.
\end{corollary}

\begin{remark} Note that the case $\a=3,\;\gamma=0,\; c>0.5,$ is covered by Theorem \ref{summing} (a).
\end{remark}

\begin{proof}
The first case readily follows from Theorem \ref{summing} (a).

Let $\a\geq 2$. If $\a>2$ or $\a=2,\;\g>0,$ for $n$ satisfying $n> 2\pi e^3+1$ and $n^{\a-2}\log^{\g}n\geq 2e^{3.5},$ invoking Proposition \ref{4} with $\xi(n):=n^{\a-2}\log^{\g}n$, we can write
\begin{align}\label{ll}
e^{n}\frac{d^n}{n^n}\cdot\frac{1}{d}<e^{n-1}\frac{d^{n-1}}{e^{\frac{1}{12}}\sqrt{2\pi}(n-1)^{n-0.5}}<\binom{d+n-2}{d-1}\leq p_{d}(n)&\leq 3e^{\frac{125n^{3-\a}}{c\log^{\g}n}}n^{2}e^{n}\frac{d^n}{n^n}\nonumber\\
&=e^{n}\frac{d^n}{n^n}e^{O(n^{3-\a}\log^{-\g}n+\log n)},
\end{align}
which gives a sharp estimate up to $e^{o(n)}$. Otherwise, when $\a=2,\;\g\leq 0,$ we obtain an extra $e^{O(n+\g n\log\log n)}$ factor at the right-hand side of \eqref{ll}, which is still $e^{O(n^{3-\a}\log^{-\g}n)}$.

Turn now to the case $\a<2$. In light of inequality \eqref{low_a}, for any $\psi=\psi(n)$ fulfilling the conditions
\begin{align*}
1\leq\psi(n)\leq \frac{n}{6},
\end{align*}
we have
\begin{align}\label{bou}
\log p_d(n+1)&\geq n\log n-\frac{n(2-\a)\log n}{\psi(n)}+\frac{n\g\log\log n}{\psi(n)}+\frac{n\log c}{\psi(n)}-\frac{3\a\log n}{2}-\frac{3\g\log \log n}{2}\nonumber\\
&-\frac{3\log C}{2}-2n\log \psi(n)+\frac{3n\log \psi(n)}{\psi(n)}-2n\log 2+\frac{2n\log 2}{\psi(n)}-\frac{3\log 2}{2}.
\end{align}
Taking $\psi(n)=\log^{\delta} n:=\log^{\max\{1,-\g\}}n$ and plugging this into \eqref{bou}, we obtain
\begin{align}\label{c}
\log p_d(n+1)&\geq n\log n-n(2-\a)\log^{1-\delta}n+\frac{n\g\log\log n}{\log^{\delta} n}+\frac{n\log c}{\log^{\delta}n}-\frac{3\a\log n}{2}-\frac{3\g\log \log n}{2}\nonumber\\
&-\frac{3\log C}{2}-2\delta n\log \log n+\frac{3\delta n\log\log n}{\log^{\delta} n}-2n\log 2+\frac{2n\log 2}{\log^{\delta} n}-\frac{3\log 2}{2},
\end{align}
which yields
\begin{align*}
\log p_d(n+1)&> n\log n+O(n\log\log n).
\end{align*}
One can see that estimate \eqref{c} is up to a constant optimal with respect to an appropriate choice of a function $\psi$. Indeed, we need to counterbalance the two main terms of \eqref{bou}, namely, $n\log n/\psi$ and $n\log \psi$. They are equal when $\psi=\log n/W(\log n)$, where $W(x)$ stands for the $W$--Lambert function, i.e. the inverse function for $ye^{y}$. The fact that $W(x)\log ^{-1}x\to 1$ as $x\to\infty$, yields $\psi(n)\sim \log n/\log\log n$ and the estimate we obtain by means of such $\psi$ is up to a constant the same as the one for $\psi(n)=\log^{\delta} n$.

At the same time according to Proposition \ref{up}, there holds
\begin{align*}
\log p_d(n)\leq n\log n+O(n).
\end{align*}
Summing up, 
$$n^ne^{O(n\log\log n)}\leq p_d(n)\leq n^ne^{O(n)}.$$
\end{proof}

\begin{acknowledgements} I am very grateful to Sergey Tikhonov for posing the problem and constantly discussing the results. I thank Miquel Saucedo for carefully reading the manuscript and giving useful comments and also the anonymous referees for their valuable suggestions.
\end{acknowledgements}

\end{document}